\documentclass[a4paper,10pt]{article}

\usepackage[english]{babel}
\usepackage{amsmath}
\usepackage{amsthm}
\usepackage{amssymb}
\usepackage{amsfonts}
\usepackage{amsxtra}
\usepackage{color}
\usepackage[breaklinks]{hyperref}
\usepackage{enumitem}
\usepackage{tocloft}

\setenumerate{topsep=0pt, itemsep=0pt}
\setitemize{topsep=0pt, itemsep=0pt}
\setdescription{topsep=0pt, itemsep=0pt}
\setlist{noitemsep}
\renewcommand{\phi}{\varphi}
\renewcommand{\epsilon}{\varepsilon}

\theoremstyle{plain}
\newtheorem{theorem}{Theorem}[section]
\newtheorem*{question}{Question}
\newtheorem*{thmA}{Theorem A}
\newtheorem*{thmB}{Theorem B}
\newtheorem*{thmC}{Theorem C}
\newtheorem*{thmD}{Theorem D}
\newtheorem*{thmE}{Theorem E}
\newtheorem*{thmF}{Theorem F}
\newtheorem{proposition}[theorem]{Proposition}
\newtheorem{lemma}[theorem]{Lemma}
\newtheorem{corollary}[theorem]{Corollary}

\theoremstyle{definition}
\newtheorem{example}[theorem]{Example}
\newtheorem{remark}[theorem]{Remark}

\theoremstyle{remark}
\newenvironment{pf}{\begin{proof}[\sc Proof]}{\end{proof}}

\newcommand{\td}{\,\mathrm{d}}
\newcommand{\Lie}{\textup{Lie}}
\newcommand{\const}{\textup{const}}

\renewcommand{\Re}{\textup{Re }}
\newcommand{\id}{\textup{id}}

\numberwithin{equation}{section}

\date{}
\title{Special functions associated to a certain fourth order differential equation}
\author{Joachim Hilgert, Toshiyuki Kobayashi\footnote{Partially supported by Grant-in-Aid for Scientific Research (B) (18340037, 22340026), Japan Society for the Promotion of Science, and the Alexander Humboldt Foundation.}, Gen Mano, Jan M\"ollers\footnote{Partially supported by the International Research Training Group 1133 ``Geometry and Analysis of Symmetries'', and the GCOE program of the University of Tokyo.}}
\makeindex

\begin{document}

\maketitle

\begin{abstract}
We develop a theory of \lq special functions\rq\ associated to a certain fourth order differential operator $\mathcal{D}_{\mu,\nu}$ on $\mathbb{R}$ depending on two parameters $\mu,\nu$. For integers $\mu,\nu\geq-1$ with $\mu+\nu\in2\mathbb{N}_0$ this operator extends to a self-adjoint operator on $L^2(\mathbb{R}_+,x^{\mu+\nu+1}\td x)$ with discrete spectrum. We find a closed formula for the generating functions of the eigenfunctions, from which we derive basic properties of the eigenfunctions such as orthogonality, completeness, $L^2$-norms, integral representations and various recurrence relations.

This fourth order differential operator $\mathcal{D}_{\mu,\nu}$ arises as the radial part of the Casimir action in the Schr\"odinger model of the minimal representation of the group $O(p,q)$, and our \lq special functions\rq\ give $K$-finite vectors.\\


\end{abstract}

\tableofcontents

\section{Introduction}

Special functions have been a powerful analytic tool in various areas of mathematics. One reason for their usefulness is the fact that there are a number of different aspects under which special functions simultaneously have good properties and admit explicit formulas (e.g. differential equations, orthonormal bases, recurrence relations, integral representations, etc.). Many of the classical special functions of this kind are associated to specific second order differential equations. In contrast, in this article we initiate the study of \lq special functions\rq\ associated to the following fourth order differential operator \index{Dmunu@$\mathcal{D}_{\mu,\nu}$|textbf}
\begin{equation}
 \mathcal{D}_{\mu,\nu} := \frac{1}{x^2}\left((\theta+\nu)(\theta+\mu+\nu)-x^2\right)\left(\theta(\theta+\mu)-x^2\right)-\frac{(\mu-\nu)(\mu+\nu+2)}{2}\label{eq:DiffOp1}
\end{equation}
with two parameters $\mu,\nu$. Here $\theta:=x\frac{\td}{\td x}$\index{$\theta$|textbf} denotes the one-dimensional Euler operator.

The differential operator $\mathcal{D}_{\mu,\nu}$ is symmetric with respect to the parameters $\mu$ and $\nu$, i.e. $\mathcal{D}_{\mu,\nu}=\mathcal{D}_{\nu,\mu}$ (see Proposition \ref{prop:DiffOpProperties}~(1)). Further, if $\mu,\nu\geq-1$ are integers of the same parity, then the operator $\mathcal{D}_{\mu,\nu}$ extends to a self-adjoint operator on $L^2(\mathbb{R}_+,x^{\mu+\nu+1}\td x)$ with discrete spectrum (see Proposition \ref{prop:DiffOpProperties}~(4)). We then ask

\begin{question}
What are the eigenfunctions of $\mathcal{D}_{\mu,\nu}$?
\end{question}

The starting point for our search for eigenfunctions of $\mathcal{D}_{\mu,\nu}$ is the fact that the principal part of $\mathcal{D}_{\mu,\nu}$ is a product of second order operators in the Weyl algebra. We will construct generating functions for eigenfunctions of $\mathcal{D}_{\mu,\nu}$ as modifications of products of generating functions for eigenfunctions of second order operators. 
More precisely, we introduce in Section \ref{sec:GenFct}
four generating functions $G_i^{\mu,\nu}(t,x)$\index{Gimunu@$G_i^{\mu,\nu}(t,x)$} ($i=1,2,3,4$). 
Of particular importance for the $L^2$-theory will be
\begin{equation*}
 G_2^{\mu,\nu}(t,x) := \frac{1}{(1-t)^{\frac{\mu+\nu+2}2}}\widetilde{I}_{\frac{\mu}{2}}\left(\frac{tx}{1-t}\right)\widetilde{K}_{\frac{\nu}{2}}\left(\frac{x}{1-t}\right),
\end{equation*}
where $\widetilde{I}_\alpha(z):=\left(\frac{z}{2}\right)^{-\alpha}I_\alpha(z)$\index{Ialpha@$\widetilde{I}_\alpha(z)$} and $\widetilde{K}_\alpha(z):=\left(\frac{z}{2}\right)^{-\alpha}K_\alpha(z)$\index{Kalpha@$\widetilde{K}_\alpha(z)$} denote the normalized $I$- and $K$-Bessel functions. $G_2^{\mu,\nu}(t,x)$ is real analytic at $t=0$ and we can define a family of real analytic functions $\Lambda_{2,j}^{\mu,\nu}(x)$ on $\mathbb{R}_+$ ($j\in\mathbb{N}_0$) as the coefficients of the Taylor expansion
\begin{equation*}
 G_2^{\mu,\nu}(t,x) = \sum_{j=0}^\infty{t^j\Lambda_{2,j}^{\mu,\nu}(x)}.
\end{equation*}
In the $L^2$-theory for the operator $\mathcal{D}_{\mu,\nu}$ we have to assume the following integrality condition
\begin{equation}
 \mu\geq\nu\geq-1\mbox{ are integers of the same parity, not both equal to $-1$.}\tag{IC1}\label{IntCond}
\end{equation}
Then we obtain the following result (for a proof see Proposition \ref{prop:DiffOpProperties}~(4), Corollary \ref{cor:Nonzeroness} and Theorem \ref{thm:EigFct}):

\begin{thmA}[$L^2$-spectrum for $\mathcal{D}_{\mu,\nu}$]
Suppose $\mu$ and $\nu$ satisfy \eqref{IntCond}. Then the differential equation
\begin{equation*}
 \mathcal{D}_{\mu,\nu}u=\lambda u
\end{equation*}
has a non-trivial solution in $L^2(\mathbb{R}_+,x^{\mu+\nu+1}\td x)$ if and only if $\lambda$ is of the form \index{$\lambda_j^{\mu,\nu}$|textbf}
\begin{equation}
 \lambda_j^{\mu,\nu} := (2j+\mu+1)^2-\frac{(\mu+1)^2}{2}-\frac{(\nu+1)^2}{2}\label{eq:EigValues}
\end{equation}
for some $j\in\mathbb{N}_0$. In this case the solution space in $L^2(\mathbb{R}_+,x^{\mu+\nu+1}\td x)$ is one-dimensional and spanned by $\Lambda_{2,j}^{\mu,\nu}$.
\end{thmA}

Further we prove the orthogonality
and closed formulas of the $L^2$-norms of 
the eigenfunctions $\Lambda_{2,j}^{\mu,\nu}(x)$ ($j\in\mathbb{N}_0$) in Corollary~\ref{cor:Norms}.
We also prove the completeness of these $L^2$-eigenfunctions under the integrality condition
\eqref{IntCond} in Corollary \ref{cor:completeness}.
Summarizing, we get:

\begin{thmB}[Orthonormal basis]
 If $\mu+\nu,\mu-\nu>-2$, then the sequence $(\Lambda_{2,j}^{\mu,\nu})_{j\in\mathbb{N}_0}$ is orthogonal in $L^2(\mathbb{R}_+,x^{\mu+\nu+1}\td x)$ with $L^2$-norms given by
 \begin{equation*}
  \|\Lambda_{2,j}^{\mu,\nu}\|_{L^2(\mathbb{R}_+,x^{\mu+\nu+1}\td x)}^2 = \frac{2^{\mu+\nu-1}\Gamma(\frac{\mu+\nu+2}{2}+j)\Gamma(\frac{\mu-\nu+2}{2}+j)}{j!(2j+\mu+1)\Gamma(j+\mu+1)}.
 \end{equation*}
 If further $\mu$ and $\nu$ satisfy \eqref{IntCond}, then $(\Lambda_{2,j}^{\mu,\nu})_{j\in\mathbb{N}_0}$ is an orthogonal basis of $L^2(\mathbb{R}_+,x^{\mu+\nu+1}\td x)$.
\end{thmB}

We now consider general eigenfunctions of
the operator $\mathcal{D}_{\mu,\nu}$
without requesting that they are square integrable. The other three generating functions $G_i^{\mu,\nu}(t,x)$ ($i=1,3,4$) are defined by similar formulas as $G_2^{\mu,\nu}(t,x)$ (see \eqref{eq:G1} - \eqref{eq:G4}), but we have to impose an additional condition for $i=3,4$:
\begin{equation}
 \mbox{$\mu$ is an odd integer $\geq1$ for $i=3,4$.}\tag{IC2}\label{IntCond2}
\end{equation}
Then again we can define families of real analytic functions $\Lambda_{i,j}^{\mu,\nu}(x)$\index{$\Lambda_{i,j}^{\mu,\nu}(x)$} on $\mathbb{R}_+$ ($j\in\mathbb{Z}$) via the Taylor ($i=1,2$), respectively Laurent ($i=3,4$) expansions of $G_i^{\mu,\nu}(t,x)$. (It is convenient to set $\Lambda_{i,j}^{\mu,\nu}(x)\equiv0$ for $i=1,2$ and $j<0$, and for $i=3,4$ and $j<-\mu$.)

The differential equation $\mathcal{D}_{\mu,\nu}u=\lambda u$ has a regular singularity at $x=0$ with characteristic exponents $\{0,-\nu,-\mu,-\mu-\nu\}$. This allows us to prove the following theorem (cf. Theorems \ref{lem:Asymptotics} and \ref{thm:EigFct}):

\begin{thmC}[Differential equation]
Suppose $\mu$ is a positive odd integer and $\nu>0$ such that $\mu-\nu\notin2\mathbb{Z}$. Then the four functions $\Lambda_{i,j}^{\mu,\nu}(x)$ ($i=1,2,3,4$) form a fundamental system of the following differential equation on $\mathbb{R}_+$:
\begin{equation*}
 \mathcal{D}_{\mu,\nu}u=\lambda_j^{\mu,\nu}u.
\end{equation*}
Furthermore, $\Lambda_{i,j}^{\mu,\nu}(x)$ has the asymptotic behavior
\begin{equation*}
 \sim c_1,c_2x^{-\nu},c_3x^{-\mu},c_4x^{-\mu-\nu}
\end{equation*}
as $x\rightarrow0$, where the closed form of the constants $c_i$ for $i=1,2,3,4$ is given in terms of Gamma functions.
\end{thmC}

The local monodromy at $x=0$ is given in Proposition \ref{prop:ParityLambda}.

The differential operator $\mathcal{D}_{\mu,\nu}$ has an irregular singularity at $x=\infty$. In Theorem \ref{lem:Asymptotics} we also describe the asymptotic behavior of the eigenfunctions $\Lambda_{i,j}^{\mu,\nu}(x)$ as $x\rightarrow\infty$. In view of Remark \ref{rem:asym-infty} this also provides the asymptotics for $x\to -\infty$.

The functions $\Lambda_{i,j}^{\mu,\nu}$ have various other remarkable properties. First we mention the integral representations in terms of Laguerre polynomials $L_n^\alpha$,\index{Lnalpha@$L_n^\alpha(x)$} which can be used to find the special values of $\Lambda_{i,j}^{\mu,\nu}(x)$ ($i=1,2$) at $\nu=\pm1$ in terms of Laguerre polynomials (see Theorem \ref{prop:IntFormulae}, Corollary \ref{cor:SpecialValue}, and Remark \ref{rem:SpecialValue2}).

\begin{thmD}[Integral representations]
 For $j\in\mathbb{N}_0$, $\Re\mu,\Re\nu>-1$ we have the following double integral representations
 \begin{align*}
  \Lambda_{1,j}^{\mu,\nu}(x) &= c_{1,j}^{\mu,\nu}\int_0^\pi{\int_0^\pi{e^{-x\cos\phi}L_j^{\frac{\mu+\nu}{2}}(x(\cos\theta+\cos\phi))\sin^\mu\theta\sin^\nu\phi\td\phi}\td\theta},\\
  \Lambda_{2,j}^{\mu,\nu}(x) &= c_{2,j}^{\mu,\nu}\int_0^\pi{\int_0^\infty{e^{-x\cosh\phi}L_j^{\frac{\mu+\nu}{2}}(x(\cos\theta+\cosh\phi))\sin^\mu\theta\sinh^\nu\phi\td\phi}\td\theta}
 \end{align*}
 with constants $c_{1,j}^{\mu,\nu}$\index{c1jmunu@$c_{1,j}^{\mu,\nu}$} and $c_{2,j}^{\mu,\nu}$\index{c2jmunu@$c_{2,j}^{\mu,\nu}$}.
\end{thmD}

Another property of the functions $\Lambda_{2,j}^{\mu,\nu}$ is that they are eigenfunctions of Meijer's $G$-transform $\mathcal{T}_{\mu,\nu}$\index{Tmunu@$\mathcal{T}_{\mu,\nu}$}. Here $\mathcal{T}_{\mu,\nu}$ is the integral operator given by
\begin{equation*}
 \mathcal{T}_{\mu,\nu}f(x) := \frac{1}{2^{\mu+\nu+1}}\int_0^\infty{G_{\mu,\nu}\left(\left(\frac{xy}{4}\right)^2\right)f(y)y^{\mu+\nu+1}\td y}
\end{equation*}
where \index{Gmunu@$G_{\mu,\nu}(t)$}\index{G2004@$G^{20}_{04}(t\vert b_1,b_2,b_3,b_4)$}
\begin{equation*}
 G_{\mu,\nu}(t) = G^{20}_{04}\left(t\left|0,-\frac{\nu}{2},-\frac{\mu}{2},-\frac{\mu+\nu}{2}\right.\right)
\end{equation*}
denotes Meijer's $G$-function (see \cite{Fox61}). It is proved in \cite{Fox61} that, if $\mu,\nu\geq-1$ are not both equal to $-1$, $\mathcal{T}_{\mu,\nu}$ is a unitary operator on $L^2(\mathbb{R}_+,x^{\mu+\nu+1}\td x)$ and $\mathcal{T}_{\mu,\nu}^2=\id$. In particular, we have the following orthogonal decomposition into $+1$ and $-1$ eigenspaces of $\mathcal{T}_{\mu,\nu}$:\index{Hmunu@$\mathcal{H}_\pm^{\mu,\nu}$|textbf}
\begin{equation*}
 L^2(\mathbb{R}_+,x^{\mu+\nu+1}\td x)=\mathcal{H}_+^{\mu,\nu}\oplus\mathcal{H}_-^{\mu,\nu}
\end{equation*}
In Section \ref{sec:GTrafo} we prove:

\begin{thmE}[Meijer's $G$-transform]
Suppose $\mu$ and $\nu$ satisfy \eqref{IntCond}. Then $\Lambda_{2,j}^{\mu,\nu}$ is an eigenfunction of Meijer's $G$-transform $\mathcal{T}_{\mu,\nu}$ for the eigenvalue $(-1)^j$. Thus, $\{\Lambda_{2,j}^{\mu,\nu}:j\in2\mathbb{N}_0\}$ forms an orthogonal basis for the Hilbert space $\mathcal{H}_+^{\mu,\nu}$ and $\{\Lambda_{2,j}^{\mu,\nu}:j\in2\mathbb{N}_0+1\}$ forms an orthogonal basis for $\mathcal{H}_-^{\mu,\nu}$.
\end{thmE}

Finally, the functions $\Lambda_{i,j}^{\mu,\nu}$ are subject to the following three recurrence relations (see Propositions \ref{prop:RecRelH}, \ref{prop:RecRelxSq}, \ref{prop:Formula1} for more details):

\begin{thmF}[Recurrence relations]
 Let $i=1,2,3,4$, $j\in\mathbb{Z}$.
 \begin{enumerate}
  \item[\textup{(1)}] The three-term recurrence relation for $x\frac{d}{dx}\Lambda_{i,j}^{\mu,\nu}$:
  \begin{multline*}
   (2j+\mu+1)\theta\Lambda_{i,j}^{\mu,\nu} = (j+1)(j+\mu+1)\Lambda_{i,j+1}^{\mu,\nu}\\
   -(2j+\mu+1)\left(\frac{\mu+\nu+2}{2}\right)\Lambda_{i,j}^{\mu,\nu}-\left(j+\frac{\mu+\nu}{2}\right)\left(j+\frac{\mu-\nu}{2}\right)\Lambda_{i,j-1}^{\mu,\nu}.
  \end{multline*}
  \item[\textup{(2)}] The recurrence relations in $\mu$ and $\nu$ (see \eqref{eq:DelEps} for the definition of the signatures $\delta(i),\epsilon(i)\in\{-1,1\}$\index{$\delta(i)$}\index{$\epsilon(i)$}):
  \begin{align*}
   \mu\left(\Lambda_{i,j}^{\mu,\nu}(x)-\Lambda_{i,j-1}^{\mu,\nu}(x)\right) &= 2\delta(i)\left(\Lambda_{i,j}^{\mu-2,\nu}(x)-\left(\frac{x}{2}\right)^2\Lambda_{i,j-2}^{\mu+2,\nu}(x)\right),\\
   \nu\left(\Lambda_{i,j}^{\mu,\nu}(x)-\Lambda_{i,j-1}^{\mu,\nu}(x)\right) &= 2\epsilon(i)\left(\Lambda_{i,j}^{\mu,\nu-2}(x)-\left(\frac{x}{2}\right)^2\Lambda_{i,j}^{\mu,\nu+2}(x)\right),\\
   \frac{\td}{\td x}\left(\Lambda_{i,j}^{\mu,\nu}(x)-\Lambda_{i,j-1}^{\mu,\nu}(x)\right) &= \delta(i)\frac{x}{2}\Lambda_{i,j-2}^{\mu+2,\nu}(x) + \epsilon(i)\frac{x}{2}\Lambda_{i,j}^{\mu,\nu+2}(x).
  \end{align*}
  \item[\textup{(3)}] The five-term recurrence relation for $x^2\Lambda_{i,j}^{\mu,\nu}(x)$:
   \begin{equation*}
    \displaystyle x^2\Lambda_{i,j}^{\mu,\nu}\in\textup{span}\{\Lambda_{i,k}^{\mu,\nu}:k=j-2,\ldots,j+2\}
   \end{equation*}
   if $j\neq-\frac{\mu-1}{2},-\frac{\mu+1}{2},-\frac{\mu+3}{2}$.
 \end{enumerate}
\end{thmF}

Since the generating functions $G_i^{\mu,\nu}(t,x)$ are given in terms of Bessel functions, Theorems B and E give rise to identities for the $I$- and $K$-Bessel functions and Meijer's $G$-function, some of which we could not trace in the literature (see Corollaries \ref{cor:BesselCompleteness}, \ref{cor:BesselIntFormula} and \ref{cor:BesselIntFormula2}).\\

Most of our proofs are of a purely analytic nature. A special role is played by generating functions. The approach to use the factorization of the principal part of the differential operator to construct generating functions might be useful also for other differential operators of a similar type.

Our proofs for Theorems A, B and E, however, partly rely on representation theory, and for this reason we had to assume the integrality condition \eqref{IntCond}. In fact, our operator $\mathcal{D}_{\mu,\nu}$ is induced from the Casimir action in the Schr\"odinger model for the minimal representation of the semisimple Lie group $G=O(p,q)$\index{G@$G$}\index{O(p,q)@$O(p,q)$} where $p=\mu+3$ and $q=\nu+3$ when (IC1) is satisfied. This representation is realized on the Hilbert space $L^2(C)$\index{L2C@$L^2(C)$}, where $C$\index{C@$C$} is an isotropic cone in $\mathbb{R}^{p+q-2}$ and is parameterized by bipolar coordinates
\begin{equation*}
 \mathbb{R}_+\times\mathbb{S}^{\mu+1}\times\mathbb{S}^{\nu+1}\stackrel{\sim}{\longrightarrow}C,(r,\omega,\eta)\mapsto(r\omega,r\eta),
\end{equation*}
where $\mathbb{S}^{n-1}$\index{Sn1@$\mathbb{S}^{n-1}$} denotes the unit sphere in $\mathbb{R}^n$. Hence $\mu+1$ and $\nu+1$ play the role of dimensions of certain spheres if the integrality condition \eqref{IntCond} is satisfied. In this case the functions $\Lambda_{2,j}^{\mu,\nu}$ are $K$-finite vectors of the representation and the recurrence relations in Theorem E can be explained by the Lie algebra action. The connection to Meijer's $G$-transform also arises from representation theory. In fact, the unitary integral transform $\mathcal{T}_{\mu,\nu}$ is the radial part of the unitary inversion operator $\mathcal{F}_C$ on $L^2(C)$ defined by the special value of the minimal representation at the longest Weyl group element. This operator plays the role of a \lq Fourier transform\rq\ on $C$.

The special values of our parameters at $\nu=\pm1$ are of particular interest. In this case, the functions $\Lambda_{2,j}^{\mu,\nu}(x)$ reduce essentially to Laguerre polynomials $L_j^\mu(x)$ (see Corollary \ref{cor:SpecialValue} and Remark \ref{rem:SpecialValue2}) and the fourth order differential operator $\mathcal{D}_{\mu,\nu}$
is of the form
\begin{equation*}
 \mathcal{D}_{\mu,\pm1}=\mathcal{S}_ {\mu,\pm1}^2+\const,
\end{equation*}
where $\mathcal{S}_{\mu,\pm1}$ is basically the standard second order Laguerre operator of having $\Lambda_{2,j}^{\mu,\pm1}$ as eigenfunctions (see Remark \ref{rem:DiffEqSpecialCase}).

In the case $\nu=-1$ the isotropic cone is the light cone in the Minkowski space which splits into the forward and backward cone. Correspondingly, the minimal representation of $O(p,2)$ splits into a highest weight and a lowest weight representation which have been well-studied also in physics. For example, it may be interpreted as the bound states of the hydrogen atom when $(\mu,\nu)=(1,-1)$, so that $G=O(4,2)$.

On the other hand, in the case $\nu=1$, it is noteworthy that Laguerre polynomials again occur as $K$-finite vectors of the minimal representations of the groups $G=O(p,4)$ with even integers $p\geq4$. In this case, the minimal representation was discovered not so long ago (see \cite{Kos90} for $SO(4,4)$), and it is known to admit neither highest nor lowest weight vectors.\\

Notation: $\mathbb{N}_0=\{0,1,2,\ldots\}$, $\mathbb{Z}=\{\ldots,-2,-1,0,1,2,\ldots\}$, $\mathbb{R}_+=\{x\in\mathbb{R}:x>0\}$.

\section{The fourth order differential operator $\mathcal{D}_{\mu,\nu}$}\label{sec:DiffOp}

In this section we collect basic properties on the fourth order differential operator $\mathcal{D}_{\mu,\nu}$\index{Dmunu@$\mathcal{D}_{\mu,\nu}$} with two parameters $\mu,\nu\in\mathbb{C}$ which was introduced in \eqref{eq:DiffOp1}.

\begin{proposition}\label{prop:DiffOpProperties}
\begin{enumerate}
 \item[\textup{(1)}] $\mathcal{D}_{\mu,\nu}=\mathcal{D}_{\nu,\mu}$.
 \item[\textup{(2)}] $\mathcal{D}_{\mu,\nu}u=\lambda u$ is a differential equation with regular singularity at $x=0$. The characteristic exponents are $0,-\mu,-\nu,-\mu-\nu$.
 \item[\textup{(3)}] If $\mu,\nu\in\mathbb{R}$, then $\mathcal{D}_{\mu,\nu}$ is a symmetric unbounded operator on the Hilbert space $L^2(\mathbb{R}_+,x^{\mu+\nu+1}\td x)$. If, in addition, $\mu+\nu,\mu-\nu>-2$, then the strictly increasing sequence $(\lambda_j^{\mu,\nu})_{j\in\mathbb{N}_0}$\index{$\lambda_j^{\mu,\nu}$} defined in \eqref{eq:EigValues} belongs to the $L^2$-spectrum of $\mathcal{D}_{\mu,\nu}$.
 \item[\textup{(4)}] If $\mu$ and $\nu$ satisfy the integrality condition \eqref{IntCond}, then $\mathcal{D}_{\mu,\nu}$ extends to a self-adjoint operator on $L^2(\mathbb{R}_+,x^{\mu+\nu+1}\td x)$. It has discrete spectrum given precisely by $(\lambda_j^{\mu,\nu})_{j\in\mathbb{N}_0}$. Furthermore, every $L^2$-eigenspace is one-dimensional.
 \item[\textup{(5)}] In the special cases where $\nu=\pm1$ the differential operator $\mathcal{D}_{\mu,\nu}$ collapses to
  \begin{equation*}
   \mathcal{D}_{\mu,\pm1} = \mathcal{S}_{\mu,\pm1}^2 - C_{\mu,\pm1},
  \end{equation*}
  where \index{Smupm1@$\mathcal{S}_{\mu,\pm1}$|textbf}\index{Cmupm1@$C_{\mu,\pm1}$|textbf}
  \begin{align*}
   \mathcal{S}_{\mu,-1} &:= \frac{1}{x}\left(\theta(\theta+\mu)-x^2\right), & C_{\mu,-1} &:= \frac{(\mu+1)^2}{2},\\
   \mathcal{S}_{\mu,+1} &:= \frac{1}{x}\left(\theta(\theta+\mu+2)+\mu+1-x^2\right), & C_{\mu,+1} &:= \frac{(\mu+1)^2}{2}+2.
  \end{align*}
\end{enumerate}
\end{proposition}

\begin{pf}
\begin{description}
\item[\normalfont\textup{(1)}] A simple computation shows that
 \begin{multline}
  \mathcal{D}_{\mu,\nu} = \frac{1}{x^2}\theta(\theta+\mu)(\theta+\nu)(\theta+\mu+\nu)+x^2\\
  -2\left(\theta^2+(\mu+\nu+2)\theta+\frac{(\mu+\nu+2)(\mu+\nu+4)}{4}\right),\label{eq:DiffOp2}
 \end{multline}
 whence $\mathcal{D}_{\mu,\nu}=\mathcal{D}_{\nu,\mu}$.
\item[\normalfont\textup{(2)}] It follows from \eqref{eq:DiffOp2} that
 \begin{equation*}
  x^2(\mathcal{D}_{\mu,\nu}-\lambda)\equiv\theta(\theta+\mu)(\theta+\nu)(\theta+\mu+\nu)\ \ \ \ \ (\textup{mod }x\cdot\mathbb{C}[x,\theta]),
 \end{equation*}
 where $\mathbb{C}[x,\theta]$ denotes the left $\mathbb{C}[x]$-module generated by $1,\theta,\theta^2,\ldots$ in the Weyl algebra $\mathbb{C}[x,\frac{\td}{\td x}]$. Therefore, the differential equation $\mathcal{D}_{\mu,\nu}u=\lambda u$ has a regular singularity at $x=0$, and its characteristic equation is given by
 \begin{equation*}
  s(s+\mu)(s+\nu)(s+\mu+\nu)=0.
 \end{equation*}
 Hence the second statement is proved.
\item[\normalfont\textup{(3)}] The formal adjoint of $\theta$ on $L^2(\mathbb{R}_+,x^{\mu+\nu+1}\td x)$ is given by
 \begin{equation*}
  \theta^*=-\theta-(\mu+\nu+2).
 \end{equation*}
It is easily seen from the expression \eqref{eq:DiffOp2} that $\mathcal{D}_{\mu,\nu}$ is a symmetric operator on the same Hilbert space. We shall postpone the proof of the assertion that $\lambda_j^{\mu,\nu}$ is an $L^2$-eigenvalue of $\mathcal{D}_{\mu,\nu}$ until Corollary \ref{prop:L2functions} and Theorem \ref{thm:EigFct}, where we give explicit eigenfunctions.
\item[\normalfont\textup{(4)}] The proof of this statement is based on representation theory. We will postpone it to Section \ref{sec:RepTh}, where we explain how the fourth order differential operator $\mathcal D_{\mu,\nu}$ arises from the Schr\"odinger model (the $L^2$-model) of the minimal representation of the indefinite orthogonal group.
\item[\normalfont\textup{(5)}] For $\nu=-1$ this follows easily from \eqref{eq:DiffOp1} by the commutator relation $[\theta,\frac{1}{x}]=-\frac{1}{x}$. The case $\nu=+1$ is also obtained by a simple computation.\qedhere
\end{description}
\end{pf}

\begin{remark}
It is likely that $\mathcal{D}_{\mu,\nu}$ still extends to a self-adjoint operator on the same
Hilbert space without assuming the integrality condition \eqref{IntCond}. For example, for $\nu=\pm1$ and arbitrary $\mu>-1$ we will construct $L^2$-eigenfunctions $\Lambda_{2,j}^{\mu,\pm1}$ of $\mathcal{D}_{\mu,\pm1}$ which are basically Laguerre polynomials (see Corollary \ref{cor:SpecialValue} and Remark \ref{rem:SpecialValue2}). Hence, they form a basis of the corresponding $L^2$-space and it follows that $\mathcal{D}_{\mu,\pm1}$ is self-adjoint with discrete spectrum. However, our proof for $\nu\not=\pm 1$ uses unitary representation theory of the semisimple Lie group $O(\mu+3,\nu+3)$ and hence uses condition \eqref{IntCond} in a crucial way.
\end{remark}

\section{The generating functions $G_i^{\mu,\nu}(t,x)$}\label{sec:GenFct}

In order to determine eigenfunctions of the $\mathcal{D}_{\mu,\nu}$ we define the following generating functions $G_i^{\mu,\nu}(t,x)$\index{Gimunu@$G_i^{\mu,\nu}(t,x)$|textbf} $i=1,2,3,4$
\begin{align}
 G_1^{\mu,\nu}(t,x) &:= \frac{1}{(1-t)^{\frac{\mu+\nu+2}2}}\widetilde{I}_{\frac{\mu}{2}}\left(\frac{tx}{1-t}\right)\widetilde{I}_{\frac{\nu}{2}}\left(\frac{x}{1-t}\right),\label{eq:G1}\\
 G_2^{\mu,\nu}(t,x) &:= \frac{1}{(1-t)^{\frac{\mu+\nu+2}2}}\widetilde{I}_{\frac{\mu}{2}}\left(\frac{tx}{1-t}\right)\widetilde{K}_{\frac{\nu}{2}}\left(\frac{x}{1-t}\right),\label{eq:G2}\\
 G_3^{\mu,\nu}(t,x) &:= \frac{1}{(1-t)^{\frac{\mu+\nu+2}2}}\widetilde{K}_{\frac{\mu}{2}}\left(\frac{tx}{1-t}\right)\widetilde{I}_{\frac{\nu}{2}}\left(\frac{x}{1-t}\right),\label{eq:G3}\\
 G_4^{\mu,\nu}(t,x) &:= \frac{1}{(1-t)^{\frac{\mu+\nu+2}2}}\widetilde{K}_{\frac{\mu}{2}}\left(\frac{tx}{1-t}\right)\widetilde{K}_{\frac{\nu}{2}}\left(\frac{x}{1-t}\right).\label{eq:G4}
\end{align}
Here \index{Ialpha@$\widetilde{I}_\alpha(z)$|textbf} \index{Kalpha@$\widetilde{K}_\alpha(z)$|textbf}
\begin{align}
 \widetilde{I}_\alpha(z) &:=\left(\frac{z}{2}\right)^{-\alpha}I_\alpha(z) = \sum_{n=0}^\infty{\frac{1}{\Gamma(n+\alpha+1)n!}\left(\frac{z}{2}\right)^{2n}},\label{eq:DefIBessel}\\
 \widetilde{K}_\alpha(z) &:=\left(\frac{z}{2}\right)^{-\alpha}K_\alpha(z) = \left(\frac{z}{2}\right)^{-\alpha}\frac{\pi}{2\sin\alpha\pi}(I_{-\alpha}(z)-I_\alpha(z))\label{eq:DefKBessel}
\end{align}
denote the normalized $I$- and $K$-Bessel functions.

Let us state the differential equations for the generating functions which we will make use of later.

\begin{lemma}[Differential equations for the generating functions]\label{lem:GenFctPDEs}
The functions $G_i^{\mu,\nu}(t,x)$, $i=1,2,3,4$, satisfy the following three differential equations:
\begin{enumerate}
 \item[\textup{(1)}] The fourth order partial differential equation
  \begin{equation*}
   \left(\mathcal{D}_{\mu,\nu}\right)_xu(t,x)=\left(4\theta_t^2+4(\mu+1)\theta_t+\frac{(\mu-\nu)(\mu+\nu+2)}{2}\right)u(t,x).
  \end{equation*}
 \item[\textup{(2)}] The second order partial differential equation
  \begin{multline*}
   (2\theta_t+\mu+1)\left(\theta_x+\frac{\mu+\nu+2}{2}\right)u(t,x)\\
   =\left(\frac{1}{t}\theta_t(\theta_t+\mu) - t\left(\theta_t+\frac{\mu+\nu+2}{2}\right)\left(\theta_t+\frac{\mu-\nu+2}{2}\right)\right)u(t,x).
  \end{multline*}
 \item[\textup{(3)}] The fifth order ordinary differential equation in $t$
  \begin{align*}
   & 8x^2\left(\theta_t+\frac{\mu-1}{2}\right)\left(\theta_t+\frac{\mu+1}{2}\right)\left(\theta_t+\frac{\mu+3}{2}\right)u(t,x)\\
   &\ \ \ \ \ =\left[\frac{2}{t^2}\theta_t\left(\theta_t-1\right)\left(\theta_t+\mu-1\right)\left(\theta_t+\mu\right)\left(\theta_t+\frac{\mu-5}{2}\right)\right.\\
   &\ \ \ \ \ \ \ \ \ \ - \frac{8}{t}\theta_t\left(\theta_t+\mu\right)\left(\theta_t+\frac{\mu-3}{2}\right)\left(\theta_t+\frac{\mu}{2}\right)\left(\theta_t+\frac{\mu+1}{2}\right)\\
   &\ \ \ \ \ \ \ \ \ \ + 2\left(\theta_t+\frac{\mu+1}{2}\right)\left(a\theta_t^4+b\theta_t^3+c\theta_t^2+d\theta_t+e\right)\\
   &\ \ \ \ \ \ \ \ \ \ - 8t\left(\theta_t+\frac{\mu+1}{2}\right)\left(\theta_t+\frac{\mu+2}{2}\right)\left(\theta_t+\frac{\mu+5}{2}\right)\\
   &\ \ \ \ \ \ \ \ \ \ \ \ \ \ \ \ \ \ \ \ \ \ \ \ \ \ \ \ \ \ \times\left(\theta_t+\frac{\mu+\nu+2}{2}\right)\left(\theta_t+\frac{\mu-\nu+2}{2}\right)\\
   &\ \ \ \ \ \ \ \ \ \ + 2t^2\left(\theta_t+\frac{\mu+7}{2}\right)\left(\theta_t+\frac{\mu+\nu+2}{2}\right)\left(\theta_t+\frac{\mu-\nu+2}{2}\right)\\
   &\ \ \ \ \ \ \ \ \ \ \ \ \ \ \ \ \ \ \ \ \ \ \ \ \ \ \ \ \ \times\left.\left(\theta_t+\frac{\mu+\nu+4}{2}\right)\left(\theta_t+\frac{\mu-\nu+4}{2}\right)\right]u(t,x),
  \end{align*}
  where we set
  \begin{align*}
   a &= 6,\\
   b &= 12(\mu+1),\\
   c &= \frac{1}{2}(17\mu^2-\nu^2+36\mu+8),\\
   d &= \frac{1}{2}(\mu+1)(5\mu^2-\nu^2+12\mu-4),\\
   e &= \frac{1}{4}(\mu-1)(\mu+2)(\mu+\nu+2)(\mu-\nu+2).
  \end{align*}
\end{enumerate}
\end{lemma}

\begin{pf}
The proof consists of straightforward verifications using the definition of $G_i^{\mu,\nu}(t,x)$ and the differential equation
\begin{equation}
 (\theta^2+2\alpha\theta-z^2)u=0\label{eq:BesselDiffEq}
\end{equation}
for the $I$- and $K$-Bessel functions $\widetilde{I}_\alpha(z)$, $\widetilde{K}_\alpha(z)$ (see \cite[Chapter III.7]{Wat44}).
\end{pf}

We will also need three recurrence relations for the functions $G_i^{\mu,\nu}(t,x)$. To state the formulas in a uniform way we put\index{$\delta(i)$|textbf}\index{$\epsilon(i)$|textbf}
\begin{align}
 \delta(i) &= \left\{\begin{array}{ll}+1&\mbox{for $i=1,2$,}\\-1&\mbox{for $i=3,4$,}\end{array}\right. & \epsilon(i) &= \left\{\begin{array}{ll}+1&\mbox{for $i=1,3$,}\\-1&\mbox{for $i=2,4$.}\end{array}\right.
 \label{eq:DelEps}
\end{align}

\begin{lemma}[Recurrence relations for the generating functions]\label{lem:GenFctRecRels}
The functions $G_i^{\mu,\nu}(t,x)$, $i=1,2,3,4$, satisfy the following three recurrence relations:
\begin{enumerate}
 \item[\textup{(1)}] The recurrence relation in $\mu$
  \begin{equation*}
   \mu(1-t)G_i^{\mu,\nu}(t,x) = 2\delta(i)\left(G_i^{\mu-2,\nu}(t,x)-\left(\frac{tx}{2}\right)^2 G_i^{\mu+2,\nu}(t,x)\right).
  \end{equation*}
 \item[\textup{(2)}] The recurrence relation in $\nu$
  \begin{equation*}
   \nu(1-t)G_i^{\mu,\nu}(t,x) = 2\epsilon(i)\left(G_i^{\mu,\nu-2}(t,x)-\left(\frac{x}{2}\right)^2 G_i^{\mu,\nu+2}(t,x)\right).
  \end{equation*}
 \item[\textup{(3)}] The recurrence relation in $\mu$ and $\nu$
  \begin{equation*}
   (1-t)\frac{\td}{\td x}G_i^{\mu,\nu}(t,x) = \delta(i)\frac{t^2x}{2}G_i^{\mu+2,\nu} + \epsilon(i)\frac{x}{2}G_i^{\mu,\nu+2}.
  \end{equation*}
\end{enumerate}
\end{lemma}

\begin{pf}
\begin{description}
 \item[\normalfont(1) and (2):] Use the recurrence relations for the $I$- and $K$-Bessel functions (see e.g. \cite[III.71~(1)]{Wat44})
  \begin{align*}
   \alpha\widetilde{I}_\alpha(x) &= \widetilde{I}_{\alpha-1}(x)-\left(\frac{x}{2}\right)^2 \widetilde{I}_{\alpha+1}(x),\\
   \alpha\widetilde{K}_\alpha(x) &= \left(\frac{x}{2}\right)^2\widetilde{K}_{\alpha+1}(x)-\widetilde{K}_{\alpha-1}(x).
  \end{align*}
 \item[\normalfont(3)] In view of the formulas (cf. \cite[III.71~(6)]{Wat44})
  \begin{align}
   \frac{\td}{\td x}\widetilde{I}_\alpha(x) &= \frac{x}{2}\widetilde{I}_{\alpha+1}(x), & \frac{\td}{\td x}\widetilde{K}_\alpha(x) &= -\frac{x}{2}\widetilde{K}_{\alpha+1}(x),\label{eq:BesselDiffFormulas}
  \end{align}
  the equation is evident.\qedhere
 \end{description}
\end{pf}

\begin{lemma}[Local monodromy of the generating functions]\label{lem:GenFctParity}
Suppose $\mu,\nu\notin2\mathbb{Z}$. Then we have the following formula for the functions $G_i^{\mu,\nu}(t,e^{-\pi\sqrt{-1}}x)$ ($i=1,2,3,4$):
\begin{equation*}
 \left(\begin{array}{c}G_1^{\mu,\nu}\\G_2^{\mu,\nu}\\G_3^{\mu,\nu}\\G_4^{\mu,\nu}\end{array}\right)(t,e^{\pi\sqrt{-1}}x) = \left(\begin{array}{cccc}1 & 0 & 0 & 0\\b_\nu & a_\nu & 0 & 0\\b_\mu & 0 & a_\mu & 0\\b_\mu b_\nu & a_\nu b_\mu & a_\mu b_\nu & a_\mu a_\nu\end{array}\right)\left(\begin{array}{c}G_1^{\mu,\nu}\\G_2^{\mu,\nu}\\G_3^{\mu,\nu}\\G_4^{\mu,\nu}\end{array}\right)(t,x)
\end{equation*}
where \index{aalpha@$a_\alpha$|textbf} \index{balpha@$b_\alpha$|textbf}
\begin{align*}
 a_\alpha &:= e^{-\alpha\pi\sqrt{-1}}, & b_\alpha &:= \frac{\Gamma(1-\frac{\alpha}{2})\Gamma(\frac{\alpha}{2})}{2}\left(e^{-\alpha\pi\sqrt{-1}}-1\right).
\end{align*}
\end{lemma}

\begin{pf}
From the definitions \eqref{eq:DefIBessel} and \eqref{eq:DefKBessel} of the Bessel functions it follows that
\begin{align*}
 \widetilde{I}_\alpha(e^{\pi\sqrt{-1}}x) &= \widetilde{I}_\alpha(x),\\
 \widetilde{K}_\alpha(e^{\pi\sqrt{-1}}x) &= a_{2\alpha}\widetilde{K}_\alpha(x)+b_{2\alpha}\widetilde{I}_\alpha(x),
\end{align*}
which proves the stated formula.
\end{pf}

\begin{remark}[Algebraic symmetries for the generating functions]
It is also easy to see that the generating functions satisfy the following algebraic symmetries
\begin{align*}
 G_i^{\mu,\nu}(t,x) &= G_i^{\nu,\mu}\left(\frac{1}{t},-x\right) & \mbox{($i=1,4$),}\\
 G_2^{\mu,\nu}(t,x) &= G_3^{\nu,\mu}\left(\frac{1}{t},-x\right).
\end{align*}
\end{remark}

\section{The eigenfunctions $\Lambda_{i,j}^{\mu,\nu}(x)$}\label{sec:EigFct}

The function $\widetilde{K}_{\frac{\mu}{2}}(\frac{tx}{1-t})$ is meromorphic near $t=0$ for a fixed $x>0$ if and only if $\mu$ is an odd integer. This explains the integrality condition \eqref{IntCond2}, which we will henceforth assume. Then the generating functions $G_i^{\mu,\nu}$ are meromorphic near $t=0$ and give rise to sequences $(\Lambda_{i,j}^{\mu,\nu})_{j\in\mathbb{Z}}$\index{$\Lambda_{i,j}^{\mu,\nu}(x)$|textbf} of functions on $\mathbb{R}_+$ as coefficients of the Laurent expansions
\begin{align}
 G_i^{\mu,\nu}(t,x) &= \sum_{j=-\infty}^\infty{t^j\Lambda_{i,j}^{\mu,\nu}(x)}, & i=1,2,3,4.\label{eq:FctDef}
\end{align}
Since $\widetilde{I}_{\frac{\mu}{2}}(z)$ is an entire function and $\widetilde{K}_{\frac{\mu}{2}}(z)$ has a pole of order $\mu$ at $z=0$ if $\mu\geq1$ is an odd integer we immediately obtain
\begin{align*}
 \Lambda_{1,j}^{\mu,\nu} = \Lambda_{2,j}^{\mu,\nu} &= 0 && \mbox{for $j<0$,}\\
 \Lambda_{3,j}^{\mu,\nu} = \Lambda_{4,j}^{\mu,\nu} &= 0 && \mbox{for $j<-\mu$.}
\end{align*}
This allows us to calculate the functions $\Lambda_{i,j}^{\mu,\nu}$ as follows:
\begin{equation}
 \Lambda_{i,j}^{\mu,\nu}(x) = \left\{\begin{array}{ll}\displaystyle\frac{1}{j!}\left.\frac{\partial^j}{\partial t^j}\right|_{t=0}G_i^{\mu,\nu}(t,x) &\mbox{if $i=1,2$, $j\geq0$,}\\\displaystyle\frac{1}{(j+\mu)!}\left.\frac{\partial^{j+\mu}}{\partial t^{j+\mu}}\right|_{t=0}t^{\mu}G_i^{\mu,\nu}(t,x) &\mbox{if $i=3,4$, $j\geq-\mu$.}\end{array}\right.\label{eq:DefAsDerivative}
\end{equation}

\begin{example}\label{ex:3Fcts}
The $L^2$-solutions $\Lambda_{2,j}^{\mu,\nu}(x)$ (i.e. the $i=2$ case) will be of special interest. Here are its first three functions.
\begin{align*}
 \Lambda_{2,0}^{\mu,\nu}(x) ={}& \frac{1}{\Gamma(\frac{\mu+2}{2})} \widetilde{K}_{\frac{\nu}{2}}(x),\\
 \Lambda_{2,1}^{\mu,\nu}(x) ={}& \frac{1}{\Gamma(\frac{\mu+2}{2})} \left(\frac{\mu+\nu+2}{2}\widetilde{K}_{\frac{\nu}{2}}(x)+\theta\widetilde{K}_{\frac{\nu}{2}}(x)\right),\\
 \Lambda_{2,2}^{\mu,\nu}(x) ={}& \frac{1}{2\Gamma(\frac{\mu+2}{2})} \left(\frac{(\mu+\nu+2)(\mu+\nu+4)}{4}\widetilde{K}_{\frac{\nu}{2}}(x)\right.\\
 &\ \ \ \ \ \ \ \ \ \ \ \ \ \ \ \ \ \ \ \ \left.+ \frac{(\mu+3)(\mu+\nu+2)}{\mu+2}\theta\widetilde{K}_{\frac{\nu}{2}}(x) + \frac{\mu+3}{\mu+2}\theta^2\widetilde{K}_{\frac{\nu}{2}}(x)\right).
\end{align*}
\end{example}

To formulate the asymptotic behavior of the functions $\Lambda_{i,j}^{\mu,\nu}(x)$ we use the Landau symbols $\mathcal{O}$ and $o$.

\begin{theorem}\label{lem:Asymptotics}
Let $\mu\in\mathbb{C}$, $\mu\neq-1,-2,-3,\ldots$ and $\nu\in\mathbb{R}$. Assume further that $j\geq0$ if $i=1,2$ and $j\geq-\mu$ if $i=3,4$.
\begin{enumerate}
\item[\textup{(1)}] The asymptotic behavior of the functions $\Lambda_{i,j}^{\mu,\nu}$ as $x\rightarrow0$ is given by
\begin{align*}
 \Lambda_{1,j}^{\mu,\nu}(x) &= \frac{(\frac{\mu+\nu+2}{2})_j}{j!\Gamma(\frac{\mu+2}{2})\Gamma(\frac{\nu+2}{2})}+o(1)\\
 \Lambda_{2,j}^{\mu,\nu}(x) &= \frac{(\frac{\mu-|\nu|+2}{2})_j}{j!\Gamma(\frac{\mu+2}{2})}\times \left\{\begin{array}{ll}\displaystyle2^{\nu-1}\Gamma\left(\frac{\nu}{2}\right)x^{-\nu}+o(x^{-\nu})&\mbox{if $\nu>0$,}\\\displaystyle-\log\left(\frac{x}{2}\right)+o\left(\log\left(\frac{x}{2}\right)\right)&\mbox{if $\nu=0$,}\\\displaystyle\frac{1}{2}\Gamma\left(-\frac{\nu}{2}\right)+o(1)&\mbox{if $\nu<0$,}\end{array}\right.\\
 \Lambda_{3,j}^{\mu,\nu}(x) &= \frac{2^{\mu-1}\Gamma(\frac{\mu}{2})(\frac{-\mu+\nu+2}{2})_{j+\mu}}{(j+\mu)!\Gamma(\frac{\nu+2}{2})}x^{-\mu}+o(x^{-\mu}),\\
 \Lambda_{4,j}^{\mu,\nu}(x) &= \frac{\Gamma(\frac{\mu}{2})(\frac{-\mu-|\nu|+2}{2})_{j+\mu}}{(j+\mu)!}\times \left\{\begin{array}{ll}\displaystyle2^{\mu+\nu-2}\Gamma\left(\frac{\nu}{2}\right)x^{-\mu-\nu}+o(x^{-\mu-\nu})&\mbox{if $\nu>0$,}\\\displaystyle-2^{\mu-1}x^{-\mu}\log\left(\frac{x}{2}\right)+o\left(x^{-\mu}\log\left(\frac{x}{2}\right)\right)&\mbox{if $\nu=0$,}\\\displaystyle2^{\mu-2}\Gamma\left(-\frac{\nu}{2}\right)x^{-\mu}+o(x^{-\mu})&\mbox{if $\nu<0$,}\end{array}\right.
\end{align*}
where $(a)_n=a(a+1)\cdots(a+n-1)$\index{$(a)_n$|textbf} is the Pochhammer symbol.
\item[\textup{(2)}] As $x\rightarrow\infty$ we have
\begin{align*}
 \Lambda_{1,j}^{\mu,\nu}(x) &= C_{1,j}^{\mu,\nu}x^{j-\frac{\nu+1}{2}}e^x\left(1+\mathcal{O}\left(\frac{1}{x}\right)\right), & \Lambda_{3,j}^{\mu,\nu}(x) &= \mathcal{O}\left(x^{j-\frac{\nu+1}{2}}e^x\right),\\
 \Lambda_{2,j}^{\mu,\nu}(x) &= C_{2,j}^{\mu,\nu}x^{j-\frac{\nu+1}{2}}e^{-x}\left(1+\mathcal{O}\left(\frac{1}{x}\right)\right), & \Lambda_{4,j}^{\mu,\nu}(x) &= \mathcal{O}\left(x^{j-\frac{\nu+1}{2}}e^{-x}\right)
\end{align*}
with constants $C_{1,j}^{\mu,\nu},C_{2,j}^{\mu,\nu}\neq0$.
\end{enumerate}
\end{theorem}

\begin{pf}
The basic ingredient for the proof is the asymptotic behavior of the Bessel functions at $x=0$ (see e.g. \cite[Chapters III and VII]{Wat44} and \cite[Chapter 4]{AAR99}):
\begin{align}
 \widetilde{I}_\alpha(0) &= \frac{1}{\Gamma(\alpha+1)},\label{eq:BesselIAsymptAt0}\\
 \widetilde{K}_\alpha(x) &= \left\{\begin{array}{ll}\frac{\Gamma(\alpha)}{2}\left(\frac{x}{2}\right)^{-2\alpha}+o(x^{-2\alpha}) &\mbox{if $\alpha>0$}\\-\log(\frac{x}{2})+o(\log(\frac{x}{2})) &\mbox{if $\alpha=0$}\\\frac{\Gamma(-\alpha)}{2}+o(1) &\mbox{if $\alpha<0$}\end{array}\right. & \mbox{as }x &\rightarrow 0,\label{eq:BesselKAsymptAt0}
\end{align}
and as $x\rightarrow\infty$
\begin{align}
\begin{split}
 \widetilde{I}_\alpha(x) &= \frac{1}{2\sqrt{\pi}}\left(\frac{x}{2}\right)^{-\alpha-\frac{1}{2}}e^x\left(1+\mathcal{O}\left(\frac{1}{x}\right)\right),\\
 \widetilde{K}_\alpha(x) &= \frac{\sqrt{\pi}}{2}\left(\frac{x}{2}\right)^{-\alpha-\frac{1}{2}}e^{-x}\left(1+\mathcal{O}\left(\frac{1}{x}\right)\right).
\end{split}\label{eq:BesselAsymptAtInfty}
\end{align}
We will also make use of the well-known expansion
\begin{equation}
 (1-t)^{-\alpha} = \sum_{j=0}^\infty{\frac{(\alpha)_j}{j!}t^j}.\label{eq:BinomExpansion}
\end{equation}
\begin{enumerate}
\item[\textup{(1)}] We show how to to calculate the asymptotic behavior at $x=0$ for the functions $\Lambda_{2,j}^{\mu,\nu}$ with $\nu>0$. The same method applies to the other cases.\\
Using the asymptotics \eqref{eq:BesselIAsymptAt0} and \eqref{eq:BesselKAsymptAt0} and the binomial expansion \eqref{eq:BinomExpansion} we find that
\begin{align*}
 \left.x^\nu G_2^{\mu,\nu}(t,x)\right|_{x=0} &= \frac{1}{(1-t)^{\frac{\mu+\nu+2}{2}}}\frac{1}{\Gamma(\frac{\mu+2}{2})}(2(1-t))^\nu\frac{\Gamma(\frac{\nu}{2})}{2}\\
 &= \sum_{j=0}^\infty{\frac{2^{\nu-1}\Gamma(\frac{\nu}{2})(\frac{\mu-\nu+2}{2})_j}{j!\Gamma(\frac{\mu+2}{2})}t^j}.
\end{align*}
In view of \eqref{eq:FctDef} this yields
\begin{equation*}
 \left.x^\nu\Lambda_{2,j}^{\mu,\nu}(x)\right|_{x=0} = \frac{2^{\nu-1}\Gamma(\frac{\nu}{2})(\frac{\mu-\nu+2}{2})_j}{j!\Gamma(\frac{\mu+2}{2})}.
\end{equation*}
\item[\textup{(2)}] Let us first treat the case $i=1,2$. With equation \eqref{eq:DefAsDerivative} it is easy to see that $\Lambda_{i,j}^{\mu,\nu}$ is a linear combination of terms of the form
\begin{equation*}
 \left\{\begin{array}{ll}\left(\theta^k\widetilde{I}_{\frac{\nu}{2}}\right)(x) & \mbox{for $i=1$,}\\\left(\theta^k\widetilde{K}_{\frac{\nu}{2}}\right)(x) & \mbox{for $i=2$}\end{array}\right.
\end{equation*}
with $0\leq k\leq j$ such that the coefficient for $k=j$ are non-zero. (In fact this can be seen in a more direct way from the recurrence relation in Proposition \ref{prop:RecRelH} and Example \ref{ex:3Fcts}.) Using \eqref{eq:BesselDiffFormulas} this simplifies to terms of the form
\begin{equation*}
 \left\{\begin{array}{ll}x^{2k}\widetilde{I}_{\frac{\nu}{2}+k}(x) & \mbox{for $i=1$,}\\x^{2k}\widetilde{K}_{\frac{\nu}{2}+k}(x) & \mbox{for $i=2$}\end{array}\right.
\end{equation*}
with $0\leq k\leq j$ and non-zero coefficient for $k=j$. Using \eqref{eq:BesselAsymptAtInfty} the leading term appears for $k=j$ and the asymptotics follow.\\

For $i=3,4$ equation \eqref{eq:DefAsDerivative} implies that $\Lambda_{i,j}^{\mu,\nu}$ is a linear combination of terms of the form
\begin{equation*}
 \left\{\begin{array}{ll}\displaystyle x^{k-\mu}\left(\theta^\ell\widetilde{I}_{\frac{\nu}{2}}\right)(x) & \mbox{for $i=3$,}\\\displaystyle x^{k-\mu}\left(\theta^\ell\widetilde{K}_{\frac{\nu}{2}}\right)(x) & \mbox{for $i=4$}\end{array}\right.
\end{equation*}
with $0\leq k+\ell\leq j+\mu$. Using \eqref{eq:BesselDiffFormulas} this simplifies to terms of the form
\begin{equation*}
 \left\{\begin{array}{ll}\displaystyle x^{k+2\ell-\mu}\widetilde{I}_{\frac{\nu}{2}+\ell}(x) & \mbox{for $i=3$,}\\\displaystyle x^{k+2\ell-\mu}\widetilde{K}_{\frac{\nu}{2}+\ell}(x) & \mbox{for $i=4$}\end{array}\right.
\end{equation*}
with $0\leq k+\ell\leq j+\mu$. Then again the claim follows from \eqref{eq:BesselAsymptAtInfty}.\qedhere
\end{enumerate}
\end{pf}

As an immediate consequence of Theorem \ref{lem:Asymptotics} we obtain

\begin{corollary}\label{prop:L2functions}
If $\mu+\nu,\mu-\nu>-2$, we have $\Lambda_{2,j}^{\mu,\nu}\in L^2(\mathbb{R}_+,x^{\mu+\nu+1}\td x)$.
\end{corollary}

From the explicit formulas for the leading terms of the functions $\Lambda_{i,j}^{\mu,\nu}$ at $x=0$ we can draw two more important corollaries.

\begin{corollary}\label{cor:Nonzeroness}
The function $\Lambda_{i,j}^{\mu,\nu}$ is non-zero if one of the following conditions is satisfied:
\begin{itemize}
\item $i=1$ and $\mu,\nu,\mu+\nu>-2$.
\item $i=2$ and $\mu+\nu,\mu-\nu>-2$.
\item $i=3,4$, $\mu$ is a positive odd integer and $\nu>-1$ such that $\mu-\nu\notin2\mathbb{Z}$.
\end{itemize}
\end{corollary}

\begin{pf}
In each case the assumption implies that the leading coeþ®¥Äient in the asymptotic expansion at $x=0$ in Theorem \ref{lem:Asymptotics} is non-zero, so that the function itself is non-zero as well.
\end{pf}

\begin{corollary}\label{cor:LinIndependence}
Suppose $\mu$ is a positive odd integer and $\nu>0$ such that $\mu-\nu\notin2\mathbb{Z}$, then for fixed $j\in\mathbb{N}_0$ the four functions $\Lambda_{i,j}^{\mu,\nu}$, $i=1,2,3,4$, are linearly independent.
\end{corollary}

\begin{pf}
The assumptions imply that the leading coefficients at $x=0$ of the functions $\Lambda_{i,j}^{\mu,\nu}(x)$ in Theorem \ref{lem:Asymptotics} never vanish and that the leading terms are distinct. Hence the asymptotic behavior near $x=0$ is different and the functions have to be linear independent.
\end{pf}

Now we can prove the main theorem of this section.

\begin{theorem}[Differential equation]\label{thm:EigFct}
For $i=1,2,3,4$, $j\in\mathbb{Z}$, the function $\Lambda_{i,j}^{\mu,\nu}$ is an eigenfunction of the fourth order differential operator $\mathcal{D}_{\mu,\nu}$ for the eigenvalue $\lambda_j^{\mu,\nu}$ (see equation \eqref{eq:EigValues}). If, in addition, $\mu$ is a positive odd integer and $\nu>0$ such that $\mu-\nu\notin2\mathbb{Z}$, then for fixed $j\in\mathbb{N}_0$ the four functions $\Lambda_{i,j}^{\mu,\nu}$, $i=1,2,3,4$, form a fundamental system of the fourth order differential equation
\begin{equation}
 \mathcal{D}_{\mu,\nu}u=\lambda_j^{\mu,\nu}u.\label{eq:DiffEq}
\end{equation}
\end{theorem}

\begin{pf}
In view of Corollary \ref{cor:LinIndependence} it only remains to show the first statement. We deduce
\begin{align}
 \mathcal{D}_{\mu,\nu}\Lambda_{i,j}^{\mu,\nu} &= \lambda_j^{\mu,\nu}\Lambda_{i,j}^{\mu,\nu} \qquad \forall j\in\mathbb{Z}\label{eq:EigFctEq}
\end{align}
from the corresponding partial differential equation for the generating function of both sides in \eqref{eq:EigFctEq}. Clearly, $(\mathcal{D}_{\mu,\nu})_xG_i^{\mu,\nu}(t,x)$ is the generating function for the left hand side of \eqref{eq:EigFctEq}. The generating function for the right hand side is calculated as follows
\begin{align*}
 & \sum_{j=-\infty}^\infty{t^j\lambda_j^{\mu,\nu}\Lambda_{i,j}^{\mu,\nu}(x)}\\
 ={}& \sum_{j=-\infty}^\infty{\left(4j^2 + 4(\mu+1)j + \frac{(\mu-\nu)(\mu+\nu+2)}{2}\right)t^j\Lambda_{i,j}^{\mu,\nu}(x)}\\
 ={}& \sum_{j=-\infty}^\infty{\left(4\theta_t^2+4(\mu+1)\theta_t + \frac{(\mu-\nu)(\mu+\nu+2)}{2}\right)t^j\Lambda_{i,j}^{\mu,\nu}(x)}\\
 ={}& \left(4\theta_t^2+4(\mu+1)\theta_t + \frac{(\mu-\nu)(\mu+\nu+2)}{2}\right)G_i^{\mu,\nu}(t,x),
\end{align*}
where $\theta_t:=t\frac{\partial}{\partial t}$. The resulting partial differential equation is
\begin{equation*}
 \left(\mathcal{D}_{\mu,\nu}\right)_xG_i^{\mu,\nu}(t,x)=\left(4\theta_t^2+4(\mu+1)\theta_t+\frac{(\mu-\nu)(\mu+\nu+2)}{2}\right)G_i^{\mu,\nu}(t,x)
\end{equation*}
which was verified in Lemma \ref{lem:GenFctPDEs}~(1).
\end{pf}

\begin{remark}
Since $\mathcal{D}_{\mu,\nu}=\mathcal{D}_{\nu,\mu}$ by Proposition \ref{prop:DiffOpProperties}~(1) and
\begin{equation*}
 \lambda_j^{\mu,\nu} = \lambda_{j+\frac{\mu-\nu}{2}}^{\nu,\mu},
\end{equation*}
Theorem \ref{thm:EigFct} implies that for $\mu-\nu\in2\mathbb{Z}$ also $\Lambda_{i,j+\frac{\mu-\nu}{2}}^{\nu,\mu}(x)$ is an eigenfunction of $\mathcal{D}_{\mu,\nu}$ for the eigenvalue $\lambda_j^{\mu,\nu}$.
\end{remark}

\begin{corollary}\label{cor:completeness}
If $\mu+\nu,\mu-\nu>-2$, then the sequence $(\Lambda_{2,j}^{\mu,\nu})_{j\in\mathbb{N}_0}$ is orthogonal in $L^2(\mathbb{R}_+,x^{\mu+\nu+1}\td x)$. If further $\mu$ and $\nu$ satisfy the integrality condition \eqref{IntCond}, then the system $(\Lambda_{2,j}^{\mu,\nu})_{j\in\mathbb{N}_0}$ forms an orthogonal basis of $L^2(\mathbb{R}_+,x^{\mu+\nu+1}\td x)$.
\end{corollary}

\begin{pf}
By Corollary \ref{prop:L2functions} the functions $\Lambda_{2,j}^{\mu,\nu}$ are contained in $L^2(\mathbb{R}_+,x^{\mu+\nu+1}\td x)$ and by Theorem \ref{thm:EigFct} each function $\Lambda_{2,j}^{\mu,\nu}$ is an eigenfunction of $\mathcal{D}_{\mu,\nu}$ for the eigenvalue $\lambda_j^{\mu,\nu}$. Further, by Proposition \ref{prop:DiffOpProperties}~(3) the operator $\mathcal{D}_{\mu,\nu}$ is symmetric on $L^2(\mathbb{R}_+,x^{\mu+\nu+1}\td x)$ and hence the orthogonality statement follows.\\
If additionally \eqref{IntCond} is satisfied, 
then by Proposition \ref{prop:DiffOpProperties}~(4) the operator $\mathcal{D}_{\mu,\nu}$ extends to a self-adjoint operator on $L^2(\mathbb{R}_+,x^{\mu+\nu+1}\td x)$ with spectrum $(\lambda_j^{\mu,\nu})_{j\in\mathbb{N}_0}$ and one-dimensional eigenspaces. This implies the completeness statement.
\end{pf}

Corollary \ref{cor:completeness} provides a completeness statement for Bessel functions we could not trace in the literature:

\begin{corollary}\label{cor:BesselCompleteness}
Fix an odd integer $N\geq1$. Then for any integer $-1\leq\nu<\frac{N}{2}$ the sequence $(\theta^j\widetilde{K}_{\frac{\nu}{2}})_{j\in\mathbb{N}_0}$ (resp. $(x^{2j}\widetilde{K}_{\frac{\nu}{2}+j})_{j\in\mathbb{N}_0}$) is a basis for $L^2(\mathbb{R}_+,x^N\td x)$. The Gram--Schmidt process applied to this sequence yields the orthogonal basis $(\Lambda_{2,j}^{N-\nu-1,\nu})_{j\in\mathbb{N}_0}$ (up to scalar factors).
\end{corollary}

\begin{pf}
Put $\mu:=N-\nu-1$, then $\mu$ and $\nu$ satisfy \eqref{IntCond}. It is an easy consequence of the definitions that $\Lambda_{2,j}^{\mu,\nu}$ can be written as a linear combination of the functions $\theta^k\widetilde{K}_{\frac{\nu}{2}}$ for $0\leq k\leq j$. (In fact this can be seen more directly from the recurrence relation in Proposition \ref{prop:RecRelH}.) Then the sequence $(\theta^j\widetilde{K}_{\frac{\nu}{2}})_j$ clearly arises from the complete sequence $(\Lambda_{2,j}^{\mu,\nu})_j$ by a base change and hence is complete. Using \eqref{eq:BesselDiffFormulas} it is also easy to see that the second series $(x^{2j}\widetilde{K}_{\frac{\nu}{2}+j})_j$ arises by a base change from the sequence $(\theta^j\widetilde{K}_{\frac{\nu}{2}})_j$. Finally, we note that both base change matrices considered are upper triangular. Thus the Gram--Schmidt process in both cases yields the orthogonal basis $\Lambda_{2,j}^{\mu,\nu}$.
\end{pf}

We end this section with a formula for the local monodromy of the functions $\Lambda_{i,j}^{\mu,\nu}(x)$ at $x=0$. This implies a parity formula with respect to $x\mapsto-x$ which can be used to determine also the asymptotic behavior as $x\rightarrow-\infty$. The monodromy formula itself is an immediate consequence of Lemma \ref{lem:GenFctParity}:

\begin{proposition}[Local monodromy at $x=0$]\label{prop:ParityLambda}
Suppose $\mu,\nu\notin2\mathbb{Z}$. Then we have the following local monodromy to the differential equation \eqref{eq:DiffEq}:
\begin{equation*}
 \left(\begin{array}{c}\Lambda_{1,j}^{\mu,\nu}\\\Lambda_{2,j}^{\mu,\nu}\\\Lambda_{3,j}^{\mu,\nu}\\\Lambda_{4,j}^{\mu,\nu}\end{array}\right)(e^{\pi\sqrt{-1}}x) = \left(\begin{array}{cccc}1 & 0 & 0 & 0\\b_\nu & a_\nu & 0 & 0\\b_\mu & 0 & a_\mu & 0\\b_\mu b_\nu & a_\nu b_\mu & a_\mu b_\nu & a_\mu a_\nu\end{array}\right)\left(\begin{array}{c}\Lambda_{1,j}^{\mu,\nu}\\\Lambda_{2,j}^{\mu,\nu}\\\Lambda_{3,j}^{\mu,\nu}\\\Lambda_{4,j}^{\mu,\nu}\end{array}\right)(x)
\end{equation*}
with non-zero coefficients $a_\alpha,b_\alpha$\index{aalpha@$a_\alpha$}\index{balpha@$b_\alpha$} as in Lemma \ref{lem:GenFctParity}.
\end{proposition}

\begin{remark}\label{rem:asym-infty}
If $\nu$ is an odd integer, the functions $\Lambda_{i,j}^{\mu,\nu}(x)$ extend holomorphically to $\mathbb{C}\backslash\{0\}$, not only to its universal covering. In this case, Proposition \ref{prop:ParityLambda} expresses $\Lambda_{i,j}^{\mu,\nu}(-x)$ as linear combination of the functions $\Lambda_{k,j}^{\mu,\nu}(x)$ ($k=1,2,3,4$). The coefficients only contain $a_\alpha$ and $b_\alpha$ with $\alpha=2n+1$ an odd integer (we assumed \eqref{IntCond2}) and simplify significantly:
\begin{align*}
 a_{2n+1} &= -1, & b_{2n+1} &= (-1)^{n+1}\pi.
\end{align*}
\end{remark}

\section{Integral representations}

In this section we show that for $i=1,2$ the functions $\Lambda_{i,j}^{\mu,\nu}(x)$ have integral representations in terms of Laguerre polynomials. Recall that for $n\in\mathbb{N}_0$ and $\alpha\in \mathbb{C}$ the Laguerre polynomial $L_n^\alpha(x)$\index{Lnalpha@$L_n^\alpha(x)$|textbf} is defined by (cf. \cite[(6.2.2)]{AAR99})
\begin{equation}
 L_n^\alpha(x) = \sum_{k=0}^n{\frac{(-1)^k}{k!}{n+\alpha\choose n-k}x^k}.\label{eq:DefLagFct}
\end{equation}

\begin{theorem}[Integral representations]\label{prop:IntFormulae}
\begin{enumerate}
\item[\textup{(1)}] For $j\in\mathbb{N}_0$, $\Re\mu,\Re\nu>-1$ we have the following double integral
representations
\begin{align}
 \Lambda_{1,j}^{\mu,\nu}(x) ={}& c_{1, j}^{\mu, \nu} \int_0^\pi{\int_0^\pi{e^{-x\cos \phi} L_j^{\frac{\mu+\nu}{2}}(x(\cos \theta+\cos \phi))\sin^\mu\theta\sin^\nu\phi\td\phi}\td\theta},\label{IntFormulaGen1}\\
 \Lambda_{2,j}^{\mu,\nu}(x) ={}& c_{2, j}^{\mu, \nu} \int_0^\pi{\int_0^\infty{e^{-x\cosh \phi} L_j^{\frac{\mu+\nu}{2}}(x(\cos \theta+\cosh \phi))\sin^\mu\theta\sinh^\nu\phi\td\phi}\td\theta},\label{IntFormulaGen2}
\end{align}
with constants $c_{1,j}^{\mu,\nu}$\index{c1jmunu@$c_{1,j}^{\mu,\nu}$|textbf} and $c_{2,j}^{\mu,\nu}$\index{c2jmunu@$c_{2,j}^{\mu,\nu}$|textbf} given by
\begin{align*}
 c_{1,j}^{\mu,\nu} &:= \frac{1}{\pi\Gamma(\frac{\mu+1}{2})\Gamma(\frac{\nu+1}{2})}\quad\text{and}\quad
 c_{2,j}^{\mu,\nu} := \frac{1}{\Gamma(\frac{\mu+1}{2})\Gamma(\frac{\nu+1}{2})}.
\end{align*}
\item[\textup{(2)}] For $\nu=-1$ and $\Re\mu>-1$ we have
\begin{align}
 \Lambda_{1,j}^{\mu,-1}(x) &= c_{1, j}^{\mu, -1} \sum_{i=0}^1{\int_0^\pi{e^{-(-1)^ix} L_j^{\frac{\mu-1}{2}}(x(\cos\theta+(-1)^i))\sin^\mu\theta\td\theta}},\label{IntFormulaSpecialCase1}\\
 \Lambda_{2,j}^{\mu,-1}(x) &= c_{2, j}^{\mu, -1} \int_0^\pi{e^{-x} L_j^{\frac{\mu-1}{2}}(x(\cos\theta+1))\sin^\mu\theta\td\theta}\label{IntFormulaSpecialCase2}
\end{align}
with constants $c_{1, j}^{\mu, -1}$ and $c_{2, j}^{\mu, -1}$ given by
\begin{align*}
 c_{1,j}^{\mu,-1} &= \frac{1}{2\pi\Gamma(\frac{\mu+1}2)}\quad\text{and}\quad
 c_{2,j}^{\mu,-1} = \frac{1}{2\Gamma(\frac{\mu+1}2)}.
\end{align*}
\end{enumerate}
\end{theorem}

\begin{pf}
We will make use of the formula \eqref{eq:DefAsDerivative} for $\Lambda_{i,j}^{\mu,\nu}$ and the generating function of the Laguerre polynomials given by (see e.g. formula (6.2.4) in
\cite{AAR99})
\begin{equation}
 \sum_{n=0}^\infty{L_n^\alpha(x)t^n} = \frac{1}{(1-t)^{\alpha+1}}e^{-\frac{tx}{1-t}}.\label{eq:LaguerreGenFct}
\end{equation}
Further, we will need the following integral representations of Bessel functions for $\Re\alpha>-\frac{1}{2}$ (cf. formulas III.71~(9) and VI.15~(5) in \cite{Wat44}) \index{Ialpha@$\widetilde{I}_\alpha(z)$}\index{Kalpha@$\widetilde{K}_\alpha(z)$}
\begin{align}
 \widetilde{I}_\alpha(x) &= \frac{1}{\sqrt{\pi}\Gamma(\alpha+\frac{1}{2})}\int_0^\pi{e^{-x\cos\theta}\sin^{2\alpha}\theta\td\theta},\label{eq:IntFormulaIBessel}\\
 \widetilde{K}_\alpha(x) &= \frac{\sqrt{\pi}}{\Gamma(\alpha+\frac{1}{2})}\int_0^\infty{e^{-x\cosh\phi}\sinh^{2\alpha}\phi\td\phi}.\label{eq:IntFormulaKBessel}
\end{align}
\begin{enumerate}
\item[\textup{(1)}] Interchanging differentiation and integration we obtain the desired integral representations for $\Lambda_{1,j}^{\mu,-1}$:
\begin{align*}
 & \pi\Gamma\left(\frac{\mu+1}{2}\right)\Gamma\left(\frac{\nu+1}{2}\right)\Lambda_{1,j}^{\mu,\nu}(x)\\
 ={}& \frac{\pi\Gamma(\frac{\mu+1}{2})\Gamma(\frac{\nu+1}{2})}{j!}\left.\frac{\partial^j}{\partial t^j}\right|_{t=0} G_1^{\mu,\nu}(t,x)\\
 ={}& \frac{1}{j!}\left.\frac{\partial^j}{\partial t^j}\right|_{t=0} \frac{1}{(1-t)^{\frac{\mu+\nu+2}{2}}} \int_0^\pi{\int_0^\pi{e^{-\frac{tx}{1-t}\cos\theta}e^{-\frac{x}{1-t}\cos\phi}\sin^\mu\theta\sin^\nu\phi\td\phi}\td\theta}\\
 ={}& \int_0^\pi{\int_0^\pi{e^{-x\cos\phi}\frac{1}{j!}\left.\frac{\partial^j}{\partial t^j}\right|_{t=0}\left[\frac{1}{(1-t)^{\frac{\mu+\nu+2}{2}}}e^{-\frac{tx}{1-t}(\cos\theta+\cos\phi)}\right]}}\\
 &\ \ \ \ \ \ \ \ \ \ \ \ \ \ \ \ \ \ \ \ \ \ \ \ \ \ \ \ \ \ \ \ \ \ \ \ \ \ \ \ \ \ \ \ \ \ \ \ \ \ \ \ \ \ \ \ \ \ \ \ \ \ \ \ \ \ \ \ \ \sin^\mu\theta\sin^\nu\phi\td\phi\td\theta\\
 ={}& \int_0^\pi{\int_0^\pi{e^{-x\cos\phi} L_j^{\frac{\mu+\nu}{2}}(x(\cos\theta+\cos\phi))\sin^\mu\theta\sin^\nu\phi\td\phi}\td\theta}.
\end{align*}
For the functions $\Lambda_{2,j}^{\mu,\nu}$ we do a similar calculation:
\begin{align*}
 & \Gamma\left(\frac{\mu+1}{2}\right)\Gamma\left(\frac{\nu+1}{2}\right)\Lambda_{2,j}^{\mu,\nu}(x)\\
 ={}& \frac{\Gamma(\frac{\mu+1}{2})\Gamma(\frac{\nu+1}{2})}{j!}\left.\frac{\partial^j}{\partial t^j}\right|_{t=0} G_2^{\mu,\nu}(t,x)\\
 ={}& \frac{1}{j!} \left.\frac{\partial^j}{\partial t^j}\right|_{t=0} \frac{1}{(1-t)^{\frac{\mu+\nu+2}{2}}} \int_0^\pi{\int_0^\infty{e^{-\frac{tx}{1-t}\cos\theta}e^{-\frac{x}{1-t}\cosh\phi}}}\\
  &\ \ \ \ \ \ \ \ \ \ \ \ \ \ \ \ \ \ \ \ \ \ \ \ \ \ \ \ \ \ \ \ \ \ \ \ \ \ \ \ \ \ \ \ \ \ \ \ \ \ \ \ \ \ \ \ \ \ \ \ \ \ \ \ \ \ \sin^\mu\theta\sinh^\nu\phi\td\phi\td\theta\\
 ={}& \int_0^\pi{\int_0^\infty{e^{-x\cosh\phi}\frac{1}{j!}\left.\frac{\partial^j}{\partial t^j}\right|_{t=0}\left[\frac{1}{(1-t)^{\frac{\mu+\nu+2}{2}}}e^{-\frac{tx}{1-t}(\cos\theta+\cosh\phi)}\right]}}\\
 &\ \ \ \ \ \ \ \ \ \ \ \ \ \ \ \ \ \ \ \ \ \ \ \ \ \ \ \ \ \ \ \ \ \ \ \ \ \ \ \ \ \ \ \ \ \ \ \ \ \ \ \ \ \ \ \ \ \ \ \ \ \ \ \ \ \ \sin^\mu\theta\sinh^\nu\phi\td\phi\td\theta\\
 ={}& \int_0^\pi{\int_0^\infty{e^{-x\cosh\phi} L_j^{\frac{\mu+\nu}{2}}(x(\cos\theta+\cosh\phi))\sin^\mu\theta\sinh^\nu\phi\td\phi}\td\theta}.
\end{align*}

\item[\textup{(2)}] Using
\begin{align}
 \widetilde{I}_{-\frac{1}{2}}(x) &= \frac{1}{\sqrt{\pi}}\cosh x\quad\text{and}\quad
 \widetilde{K}_{-\frac{1}{2}}(x) = \frac{\sqrt{\pi}}{2}e^{-x},\label{eq:SpecialCaseBesselFcts}
\end{align}
similar calculations as in (1) give the second part.\qedhere
\end{enumerate}
\end{pf}

\begin{remark}
The integral representations in Theorem \ref{prop:IntFormulae}~(2) for the special case $\nu=-1$ can also be obtained from the integral representations in part (1) for $\nu>-1$ by taking the limit $\nu\rightarrow-1$. For example, to obtain the integral representation for $\Lambda_{2,j}^{\mu,-1}$ we have to verify the limit formula
\begin{equation}
 \lim_{\nu\rightarrow-1}{\frac{1}{\Gamma(\frac{\nu+1}{2})}\int_0^\infty{e^{-x\cosh\phi}\cosh^k\phi\sinh^\nu\phi\td\phi}} = \frac{1}{2}e^{-x}\label{eq:IntLimit}
\end{equation}
for $0\leq k\leq j$. For $k=0$ the identity \eqref{eq:IntFormulaKBessel} turns the left hand side into
\begin{equation*}
 \frac{1}{\sqrt{\pi}}\lim_{\nu\rightarrow-1}{\widetilde{K}_{\frac{\nu}{2}}(x)}.
\end{equation*}
The map $\alpha\mapsto\widetilde{K}_\alpha(x)$ is continuous so \eqref{eq:IntLimit} follows from \eqref{eq:SpecialCaseBesselFcts}. For $k>0$ and $\phi\geq0$ we have
\begin{align*}
 \cosh^k\phi-\cosh^0\phi &= \cosh^k\phi-1 \leq \sinh\phi \cdot p(\sinh\phi)
\end{align*}
with some polynomial $p$. Then one has to show that
\begin{equation*}
 \lim_{\nu\rightarrow-1}{\frac{1}{\Gamma(\frac{\nu+1}{2})}\int_0^\infty{e^{-x\cosh\phi}\sinh^{\nu+\ell+1}\phi\td\phi}} = \frac{1}{2}e^{-x}.
\end{equation*}
for $\ell\geq0$. But this is easily seen using the integral representation \eqref{eq:IntFormulaKBessel} and the continuity of the map $\alpha\mapsto\widetilde{K}_\alpha(x)$.
\end{remark}

As an easy application of the integral representations we give explicit expressions for the functions $\Lambda_{i,j}^{\mu,\nu}$, $i=1,2$, in the case where $\nu=-1$.

\begin{corollary}\label{cor:SpecialValue}
For $\nu=-1$ and $\mu\in\mathbb{C}$ arbitrary we have\index{Lnalpha@$L_n^\alpha(x)$} the following identity of meromorphic functions
\begin{align}
 \Lambda_{1,j}^{\mu,-1}(x) &= \frac{2^{\mu-1}\Gamma(j+\frac{\mu+1}{2})}{\pi\Gamma(j+\mu+1)}\left(e^{-x}L_j^\mu(2x)+e^{x}L_j^\mu(-2x)\right),\label{eq:Lambda1SpecialValue}\\
 \Lambda_{2,j}^{\mu,-1}(x) &= \frac{2^{\mu-1}\Gamma(j+\frac{\mu+1}{2})}{\Gamma(j+\mu+1)} e^{-x} L_j^{\mu}(2x).\label{eq:Lambda2SpecialValue}
\end{align}
\end{corollary}

\begin{pf}
For the proof we may assume that $\Re\mu>-1$. The general case $\mu\in\mathbb{C}$ then follows by meromorphic continuation. Note the identity (cf. formula 16.6~(5) in \cite{EMOT54a})
\begin{multline*}
 \int_0^1{(1-y)^{\beta-\alpha-1}y^\alpha L_n^\alpha(xy)\td y} = \frac{\Gamma(\alpha+n+1)\Gamma(\beta-\alpha)}{\Gamma(\beta+n+1)}L_n^\beta(x),\\\Re\beta>\Re\alpha>-1.
\end{multline*}
With this integral formula the substitution $y=\frac{1}{2}(1\pm\cos\theta)$ yields
\begin{align}
 \int_0^\pi{L_j^{\frac{\mu-1}{2}}(x(\cos\theta\pm1))\sin^\mu\theta\td\theta}
 ={}& 2^\mu \int_0^1{(1-y)^{\frac{\mu-1}{2}}y^{\frac{\mu-1}{2}}L_j^{\frac{\mu-1}{2}}(\pm2x\cdot y)\td y}\notag\\
 ={}& \frac{2^\mu\Gamma(j+\frac{\mu+1}{2})\Gamma(\frac{\mu+1}{2})}{\Gamma(j+\mu+1)}L_j^\mu(\pm2x).\label{eq:LegendreIntegral2}
\end{align}
Inserting this into the integral representations \eqref{IntFormulaSpecialCase1} and \eqref{IntFormulaSpecialCase2} gives \eqref{eq:Lambda1SpecialValue} and \eqref{eq:Lambda2SpecialValue}.
\end{pf}

\begin{remark}\label{rem:SpecialValue2}
The symmetry property $K_{-\alpha}=K_\alpha$ for the $K$-Bessel functions implies that $G_2^{\mu,-1}(x)=\frac{x}{2}G_2^{\mu,1}(x)$ and hence Corollary \ref{cor:SpecialValue} also allows us to compute $\Lambda_{2,j}^{\mu,1}$ explicitly:
\begin{equation}
 \Lambda_{2,j}^{\mu,1}(x) = \frac{2}{x}\Lambda_{2,j}^{\mu,-1}(x) = \frac{2^\mu\Gamma(j+\frac{\mu+1}{2})}{\Gamma(j+\mu+1)} x^{-1}e^{-x} L_j^{\mu}(2x).\label{eq:RelLambdapm1}
\end{equation}
In \cite{HKMM09b} we show that the functions $\Lambda_{2,j}^{\mu,\nu}$ always degenerate to polynomials if $\nu$ is an odd integer. For $\nu\geq3$ these special polynomials do not appear in the standard literature. Properties for these polynomials such as differential equations, orthogonality relations, completeness, recurrence relations and integral representations will be studied thoroughly in \cite{HKMM09b}.
\end{remark}

\begin{remark}\label{rem:DiffEqSpecialCase}
Corollary \ref{cor:SpecialValue} and Remark \ref{rem:SpecialValue2} suggest a relation between the fourth order differential equation $\mathcal{D}_{\mu,\nu}u=\lambda_j^{\mu,\nu}u$ in the cases where $\nu=\pm1$ and the second order differential equation (see \cite[(6.2.8)]{AAR99})
\begin{equation}
 \left(x\frac{\td^2}{\td x^2}+(\alpha+1-x)\frac{\td}{\td x}+n\right)u = 0\label{eq:LagDiffEq}
\end{equation}
for the Laguerre polynomials $L_n^\alpha(x)$. In fact, by Proposition \ref{prop:DiffOpProperties}~(5) the fourth order differential operator $\mathcal{D}_{\mu,\pm1}$ collapses to the simpler form
\begin{align*}
 \mathcal{D}_{\mu,\pm1} &= \mathcal{S}_{\mu,\pm1}^2-C_{\mu,\pm1}
\end{align*}
with second order differential operators $\mathcal{S}_{\mu,\pm1}$\index{Smupm1@$\mathcal{S}_{\mu,\pm1}$} and constants $C_{\mu,\pm1}$\index{Cmupm1@$C_{\mu,\pm1}$} defined in Proposition \ref{prop:DiffOpProperties}~(5).

For $\mu>-1$ the operator $\mathcal{S}_{\mu,-1}$ itself is self-adjoint on $L^2(\mathbb{R}_+,x^\mu\td x)$. It has discrete spectrum given by $(-(2j+\mu+1))_{j\in\mathbb{N}_0}$ and an easy calculation involving \eqref{eq:LagDiffEq} shows that $\Lambda_{2,j}^{\mu,-1}$ is an eigenfunction of $\mathcal{S}_{\mu,-1}$ for the eigenvalue $-(2j+\mu+1)$. Applying $\mathcal{S}_{\mu,-1}$ twice yields the fourth order differential equation of Theorem \ref{thm:EigFct} for $\Lambda_{2,j}^{\mu,-1}$.

The same considerations apply for $\mathcal{S}_{\mu,+1}$ and $\Lambda_{2,j}^{\mu,+1}(x)$ since we have the relations \eqref{eq:RelLambdapm1} and
\begin{equation}
 \mathcal{S}_{\mu,+1}x^{-1} = x^{-1}\mathcal{S}_{\mu,-1}.\label{eq:RelSmupm1}
\end{equation}
\end{remark}

\begin{remark}
Similarly to the calculations in \cite{HKMM09b}, where we derive explicit formulas for the functions $\Lambda_{2,j}^{\mu,\nu}$ if $\nu$ is an odd integer, one can find the following expression for $\Lambda_{3,j}^{\mu,\nu}$:
\begin{multline*}
 \Lambda_{3,j}^{\mu,\nu}(x) = \frac{2^\nu\Gamma(\frac{\nu+1}{2})}{\Gamma(\nu+1)}x^{-\mu}e^x\sum_{i=0}^{\min(\frac{\mu-1}{2},j+\mu)}{\frac{(\mu-i-1)!\Gamma(j+\mu+\frac{\nu-\mu+2}{2})}{i!(\frac{\mu-1}{2}-i)!(j+\mu-i)!\Gamma(i+\frac{\nu-\mu+2}{2})}}\\
 \times(2x)^i{_2F_2}\left(i-j-\mu,\frac{\nu+1}{2};i+\frac{\nu-\mu+2}{2},\nu+1;2x\right),
\end{multline*}
where ${_2F_2}(a_1,a_2;b_1,b_2;z)$\index{F(a_1,a_2;b_1,b_2;z)@${_2F_2}(a_1,a_2;b_1,b_2;z)$} is the hypergeometric function. Using this one can sharpen the estimates for $\Lambda_{3,j}^{\mu,\nu}(x)$ given in Theorem \ref{lem:Asymptotics}.
\end{remark}

\section{Recurrence relations}\label{sec:RecRel}

In this section we will give three types of recurrence relations for the functions $\Lambda_{i,j}^{\mu,\nu}$. Our first recurrence relation involves the first order differential operator $\mathcal{H}_\alpha$\index{Halpha@$\mathcal{H}_\alpha$|textbf} ($\alpha\in\mathbb{C}$) on $\mathbb{R}_+$, given by
\begin{equation*}
 \mathcal{H}_\alpha := \theta+\frac{\alpha+2}{2}.
\end{equation*}
If $\alpha\in\mathbb{R}$, then $\mathcal{H}_\alpha$ is a skewsymmetric operator on $L^2(\mathbb{R}_+,x^{\alpha+1}\td x)$ since $\theta^*=-\theta-(\alpha+2)$ on $L^2(\mathbb{R}_+,x^{\alpha+1}\td x)$. This will allow us to compute the $L^2$-norms for $\Lambda_{2,j}^{\mu,\nu}$ explicitly if $\mu$ and $\nu$ satisfy \eqref{IntCond}.

\begin{proposition}\label{prop:RecRelH}
For $\mu,\nu\in\mathbb{C}$, $i=1,2,3,4$ we have the following recurrence relation in $j\in\mathbb{Z}$
\begin{multline}
 (2j+\mu+1)\mathcal{H}_{\mu+\nu}\Lambda_{i,j}^{\mu,\nu}\\
 = (j+1)(j+\mu+1)\Lambda_{i,j+1}^{\mu,\nu} - \left(j+\frac{\mu+\nu}{2}\right)\left(j+\frac{\mu-\nu}{2}\right)\Lambda_{i,j-1}^{\mu,\nu}.\label{eq:RecRelH}
\end{multline}
\end{proposition}

\begin{pf}
As in the proof of Theorem \ref{thm:EigFct} we verify \eqref{eq:RecRelH} via a partial differential equation for the generating function $G_i^{\mu,\nu}$. A short calculation shows that the recurrence relation \eqref{eq:RecRelH} is equivalent to the partial differential equation
\begin{multline*}
 (2\theta_t+\mu+1)\left(\theta_x+\frac{\mu+\nu+2}{2}\right)G_i^{\mu,\nu}\\
 = \left(\frac{1}{t}\theta_t(\theta_t+\mu) - t\left(\theta_t+\frac{\mu+\nu+2}{2}\right)\left(\theta_t+\frac{\mu-\nu+2}{2}\right)\right)G_i^{\mu,\nu},
\end{multline*}
which holds by Lemma \ref{lem:GenFctPDEs}~(2).
\end{pf}

\begin{corollary}\label{cor:Norms}
If $\mu+\nu,\mu-\nu>-2$, then
\begin{equation}
 \|\Lambda_{2,j}^{\mu,\nu}\|_{L^2(\mathbb{R}_+,x^{\mu+\nu+1}\td x)}^2 = \frac{2^{\mu+\nu-1}\Gamma(j+\frac{\mu+\nu+2}{2})\Gamma(j+\frac{\mu-\nu+2}{2})}{j!(2j+\mu+1)\Gamma(j+\mu+1)}.\label{eq:Norms}
\end{equation}
\end{corollary}

\begin{pf}
We prove this by induction on $j$. For $j=0$, in view of Example \ref{ex:3Fcts} and
\begin{multline*}
 \int_0^\infty{x^{\sigma-1}K_\mu(x)K_\nu(x)dx} = \frac{2^{\sigma-3}}{\Gamma(\sigma)}\Gamma\left(\frac{\sigma+\mu+\nu}{2}\right)\Gamma\left(\frac{\sigma+\mu-\nu}{2}\right)\\
 \times\Gamma\left(\frac{\sigma-\mu+\nu}{2}\right)\Gamma\left(\frac{\sigma-\mu-\nu}{2}\right)
\end{multline*}
(see formula 10.3~(49) in \cite{EMOT54a} for $a=y=1$ and $\Re\sigma>|\Re\mu|+|\Re\nu|$), we can calculate
\begin{align*}
 \|\Lambda_{2,0}^{\mu,\nu}\|^2 &= \int_0^\infty{|\Lambda_{2,0}^{\mu,\nu}(x)|^2x^{\mu+\nu+1}\td x}\\
 &= \frac{1}{\Gamma(\frac{\mu+2}{2})^2}\int_0^\infty{|\widetilde{K}_{\frac{\nu}{2}}(x)|^2x^{\mu+\nu+1}\td x}\\
 &= \frac{2^\nu}{\Gamma(\frac{\mu+2}{2})^2}\int_0^\infty{|K_{\frac{\nu}{2}}(x)|^2x^{\mu+1}\td x}\\
 &= \frac{2^\nu}{\Gamma(\frac{\mu+2}{2})^2} \cdot \frac{2^{\mu-1}\Gamma(\frac{\mu+2}{2})^2\Gamma(\frac{\mu+\nu+2}{2})\Gamma(\frac{\mu-\nu+2}{2})}{\Gamma(\mu+2)}\\
 &= \frac{2^{\mu+\nu-1}\Gamma(\frac{\mu+\nu+2}{2})\Gamma(\frac{\mu-\nu+2}{2})}{(\mu+1)\Gamma(\mu+1)}.
\end{align*}
For the induction step we reformulate \eqref{eq:RecRelH} as
\begin{equation*}
 \mathcal{H}_{\mu+\nu}\Lambda_{2,j}^{\mu,\nu} = \frac{(j+1)(j+\mu+1)}{2j+\mu+1}\Lambda_{2,j+1}^{\mu,\nu} - \frac{(2j+\mu+\nu)(2j+\mu-\nu)}{4(2j+\mu+1)}\Lambda_{2,j-1}^{\mu,\nu}
\end{equation*}
and use the skew-symmetry of $\mathcal{H}_{\mu+\nu}$ on $L^2(\mathbb{R}_+,x^{\mu+\nu+1}\td x)$ together with the pairwise orthogonality of the functions $\Lambda_{2,j}^{\mu,\nu}$ (cf. Corollary \ref{cor:completeness}) to calculate
\begin{align*}
 {}& \|\mathcal{H}_{\mu+\nu}\Lambda_{2,j}^{\mu,\nu}\|^2
 = \left(\mathcal{H}_{\mu+\nu}\Lambda_{2,j}^{\mu,\nu}|\mathcal{H}_{\mu+\nu}\Lambda_{2,j}^{\mu,\nu}\right)
 = - \left(\Lambda_{2,j}^{\mu,\nu}|\mathcal{H}_{\mu+\nu}^2\Lambda_{2,j}^{\mu,\nu}\right)\\
 ={}& - \left(\Lambda_{2,j}^{\mu,\nu}\left|\frac{(j+1)(j+\mu+1)}{2j+\mu+1}\mathcal{H}_{\mu+\nu}
  \Lambda_{2,j+1}^{\mu,\nu}\right.\right)\\
 &\ \ \ \ \ \ \ \ \ \ \ \ \ \ \ \ \ \ \ \ \ \ \ \ \ \ \ \ \ \ \ \ \ \ \ \ \ \ \ \ \ \ \ \ \ \ \left.-\frac{(2j+\mu+\nu)(2j+\mu-\nu)}{4(2j+\mu+1)}\mathcal{H}_{\mu+\nu}\Lambda_{2,j-1}^{\mu,\nu}\right)\\
 ={}& \left(\frac{(j+1)(j+\mu+1)}{2j+\mu+1}\cdot\frac{(2(j+1)+\mu+\nu)(2(j+1)+\mu-\nu)}{4(2(j+1)+\mu+1)}
   \right.\\
 &\ \ \ \ \ \ \ \ \ \ \ \ \ \ \ \ \ \ \ \ \ \ \
 +\left.\frac{(2j+\mu+\nu)(2j+\mu-\nu)}{4(2j+\mu+1)}\cdot\frac{j((j-1)+\mu+1)}{2(j-1)+\mu+1}\right)\|
 \Lambda_{2,j}^{\mu,\nu}\|^2.
\end{align*}
On the other hand, orthogonality and recurrence relation also yield
\begin{multline*}
 \|\mathcal{H}_{\mu+\nu}\Lambda_{2,j}^{\mu,\nu}\|^2 = \left(\frac{(j+1)(j+\mu+1)}{2j+\mu+1}\right)^2\|\Lambda_{2,j+1}^{\mu,\nu}\|^2\\
 +\left(\frac{(2j+\mu+\nu)(2j+\mu-\nu)}{4(2j+\mu+1)}\right)^2\|\Lambda_{2,j-1}^{\mu,\nu}\|^2.
\end{multline*}
Putting both equalities together completes the induction.
\end{pf}

Formula \eqref{eq:Norms} implies an integral formula for the $I$- and $K$-Bessel functions which we did not find in the literature:

\begin{corollary}\label{cor:BesselIntFormula}
Let $\frac{1}{2}<\alpha<\infty$ and $\mu+\nu,\mu-\nu>-2$. Then
\begin{multline*}
 \int_0^\infty{|\widetilde{I}_{\frac{\mu}{2}}((\alpha-1)x)|^2|\widetilde{K}_{\frac{\nu}{2}}(\alpha x)|^2 x^{\mu+\nu+1}\td x}\\
 = \frac{(2\sqrt{\alpha})^{\mu+\nu+2}}{8}\sum_{j=0}^\infty{\frac{(\alpha-1)^{2j}\Gamma(j+\frac{\mu+\nu+2}{2})\Gamma(j+\frac{\mu-\nu+2}{2})}{\alpha^{2j}j!(2j+\mu+1)\Gamma(j+\mu+1)}}.
\end{multline*}
\end{corollary}

\begin{pf} The orthogonality of the $\Lambda_{2,j}^{\mu,\nu}$ gives
\begin{align*}
 \sum_{j=0}^\infty{t^{2j}\|\Lambda_{2,j}^{\mu,\nu}\|^2} &= \bigg\|\sum_{j=0}^\infty{t^j\Lambda_{2,j}^{\mu,\nu}}\bigg\|^2 = \|G_2^{\mu,\nu}(t,-)\|^2\\
 &= \frac{1}{(1-t)^{\frac{\mu+\nu+2}{2}}}\int_0^\infty{\left|\widetilde{I}_{\frac{\mu}{2}}\left(\frac{tx}{1-t}\right)\right|^2\left|\widetilde{K}_{\frac{\nu}{2}}\left(\frac{x}{1-t}\right)\right|^2x^{\mu+\nu+1}\td x},
\end{align*}
so that for $\alpha:=\frac{1}{1-t}$ the desired formula follows from \eqref{eq:Norms}.
\end{pf}

The second type of recurrence relations expresses $x^2\Lambda_{i,j}^{\mu,\nu}$ as linear combination in $\Lambda_{i,k}^{\mu,\nu}$ for $k=j-2,\ldots,j+2$. These recurrence relations are an immediate consequence of the fifth order differential equation for the generating function $G_i^{\mu,\nu}$ given in Lemma \ref{lem:GenFctPDEs}~(3).

\begin{proposition}\label{prop:RecRelxSq}
For $\mu,\nu\in\mathbb{C}$ we have
\begin{align*}
 &8\left(j+\frac{\mu-1}{2}\right)\left(j+\frac{\mu+1}{2}\right)\left(j+\frac{\mu+3}{2}\right)x^2\Lambda_{i,j}^{\mu,\nu}(x)\notag\\
 &\ \ \ \ \ =2(j+1)(j+2)(j+\mu+1)(j+\mu+2)\left(j+\frac{\mu-1}{2}\right)\Lambda_{i,j+2}^{\mu,\nu}(x)\notag\\
 &\ \ \ \ \ \ \ \ - 8(j+1)(j+\mu+1)\left(j+\frac{\mu-1}{2}\right)\left(j+\frac{\mu+2}{2}\right)\left(j+\frac{\mu+3}{2}\right)\Lambda_{i,j+1}^{\mu,\nu}(x)\notag\\
 &\ \ \ \ \ \ \ \ + 2\left(j+\frac{\mu+1}{2}\right)(aj^4+bj^3+cj^2+dj+e)\Lambda_{i,j}^{\mu,\nu}(x)\notag\\
 &\ \ \ \ \ \ \ \ - 8\left(j+\frac{\mu-1}{2}\right)\left(j+\frac{\mu}{2}\right)\left(j+\frac{\mu+3}{2}\right)\left(j+\frac{\mu+\nu}{2}\right)\notag\\
 &\ \ \ \ \ \ \ \ \ \ \ \ \ \ \ \ \ \ \ \ \ \ \ \ \ \ \ \ \ \ \ \ \ \ \ \ \ \ \ \ \ \ \ \ \ \ \ \ \ \ \ \ \ \ \ \ \ \ \ \ \ \ \ \ \ \ \ \left(j+\frac{\mu-\nu}{2}\right)\Lambda_{i,j-1}^{\mu,\nu}(x)\notag\\
 &\ \ \ \ \ \ \ \ + 2\left(j+\frac{\mu+3}{2}\right)\left(j+\frac{\mu+\nu-2}{2}\right)\left(j+\frac{\mu-\nu-2}{2}\right)\notag\\
 &\ \ \ \ \ \ \ \ \ \ \ \ \ \ \ \ \ \ \ \ \ \ \ \ \ \ \ \ \ \ \ \ \ \ \ \ \ \ \ \ \ \ \ \ \ \ \ \ \times\left(j+\frac{\mu+\nu}{2}\right)\left(j+\frac{\mu-\nu}{2}\right)\Lambda_{i,j-2}^{\mu,\nu}(x)
\end{align*}
with $a,b,c,d,e$ as in Lemma \ref{lem:GenFctPDEs}~(3).
\end{proposition}

\begin{remark}
For $j\neq-\frac{\mu-1}{2},-\frac{\mu+1}{2},-\frac{\mu+3}{2}$ the recurrence relation of Proposition \ref{prop:RecRelxSq} can be rewritten as
\begin{equation*}
 x^2\Lambda_{i,j}^{\mu,\nu}(x) = \sum_{k=-2}^{2}{a_{i,j}^{\mu,\nu}(k)\Lambda_{i,j+k}^{\mu,\nu}(x)}
\end{equation*}
with constants $a_{i,j}^{\mu,\nu}(k)$.
\end{remark}

The last set of recurrence relations in the parameters $\mu$ and $\nu$ are again immediate with the corresponding differential equations for the generating functions which have already been stated in Lemma \ref{lem:GenFctRecRels}:

\begin{proposition}\label{prop:Formula1}
Let $\mu,\nu\in\mathbb{C}$. With $\delta(i)$\index{$\delta(i)$}, $\epsilon(i)$\index{$\epsilon(i)$} as in \eqref{eq:DelEps} we have
\begin{enumerate}
 \item[\textup{(1)}] The recurrence relation in $\mu$
  \begin{equation*}
   \mu\left(\Lambda_{i,j}^{\mu,\nu}(x)-\Lambda_{i,j-1}^{\mu,\nu}(x)\right) = 2\delta(i)\left(\Lambda_{i,j}^{\mu-2,\nu}(x)-\left(\frac{x}{2}\right)^2\Lambda_{i,j-2}^{\mu+2,\nu}(x)\right).
  \end{equation*}
 \item[\textup{(2)}] The recurrence relation in $\nu$
  \begin{equation*}
   \nu\left(\Lambda_{i,j}^{\mu,\nu}(x)-\Lambda_{i,j-1}^{\mu,\nu}(x)\right) = 2\epsilon(i)\left(\Lambda_{i,j}^{\mu,\nu-2}(x)-\left(\frac{x}{2}\right)^2\Lambda_{i,j}^{\mu,\nu+2}(x)\right).
  \end{equation*}
 \item[\textup{(3)}] The recurrence relation in $\mu$ and $\nu$
  \begin{equation*}
   \frac{\td}{\td x}\left(\Lambda_{i,j}^{\mu,\nu}(x)-\Lambda_{i,j-1}^{\mu,\nu}(x)\right) = \delta(i)\frac{x}{2}\Lambda_{i,j-2}^{\mu+2,\nu}(x) + \epsilon(i)\frac{x}{2}\Lambda_{i,j}^{\mu,\nu+2}(x).
  \end{equation*}
\end{enumerate}
\end{proposition}

\section{Meijer's $G$-transform}\label{sec:GTrafo}

The main result of this section is that the functions $\Lambda_{2,j}^{\mu,\nu}$ are eigenfunctions of a special type of Meijer's $G$-transform $\mathcal{T}_{\mu,\nu}$\index{Tmunu@$\mathcal{T}_{\mu,\nu}$|textbf} if $\mu$ and $\nu$ satisfy \eqref{IntCond}. Our proof of this result is based on representation theory and will be given in Section \ref{sec:RepTh}. We do not know whether the statement also holds for general $\mu\geq\nu\geq-1$.

For $f\in\mathcal{C}_c^\infty(\mathbb{R}_+)$ the integral transform $\mathcal{T}_{\mu,\nu}$ it is defined by
\begin{equation*}
 \mathcal{T}_{\mu,\nu}f(x) := \frac{1}{2^{\mu+\nu+1}}\int_0^\infty{G_{\mu,\nu}\left(\left(\frac{xy}{4}\right)^2\right)f(y)y^{\mu+\nu+1}\td y}.
\end{equation*}
Here
\begin{equation*}
 G_{\mu,\nu}(t):=G^{20}_{04}\left(t\left|0,-\frac{\nu}{2},-\frac{\mu}{2},-\frac{\mu+\nu}{2}\right.\right)
\end{equation*}
denotes Meijer's $G$-function, which was first systematically investigated by Fox in \cite{Fox61}.\index{Gmunu@$G_{\mu,\nu}(t)$|textbf}\index{G2004@$G^{20}_{04}(t\vert b_1,b_2,b_3,b_4)$|textbf} $G_{\mu,\nu}$ satisfies the fourth order differential equation (see \cite[5.4~(1)]{EMOT81})
\begin{equation}
 \theta\left(\theta+\frac{\mu}{2}\right)\left(\theta+\frac{\nu}{2}\right)\left(\theta+\frac{\mu+\nu}{2}\right)u(t)=tu(t).\label{eq:GDiffEq}
\end{equation}
For $\nu=-1$ the $G$-function reduces to a $J$-Bessel function\index{Jalpha@$J_\alpha(z)$|textbf}
\begin{equation*}
 G^{\mu,-1}(t) = t^{-\frac{\mu}{4}}J_\mu(4t^{\frac{1}{4}})
\end{equation*}
and the operator $\mathcal{T}_{\mu,-1}$ becomes a Hankel transform.

The following statement can also be obtained from representation theory in the case that $\mu$ and $\nu$ satisfy \eqref{IntCond} (see Section \ref{sec:RepTh} for a proof).

\begin{proposition}\label{prop:GTrafoProperties}
Suppose $\mu+\nu>-2$.
\begin{enumerate}
 \item[\textup{(1)}] $\mathcal{T}_{\mu,\nu}$ extends to a unitary involutive operator on $L^2(\mathbb{R}_+,x^{\mu+\nu+1}dx)$.
 \item[\textup{(2)}] The $G$-transform $\mathcal{T}_{\mu,\nu}$ commutes with the fourth order differential operator $\mathcal{D}_{\mu,\nu}$.
\end{enumerate}
\end{proposition}

\begin{pf}
\begin{enumerate}
 \item[\textup{(1)}] It is proved in \cite{Fox61} that
 \begin{equation*}
  Tf(r) := \frac{1}{c}\int_0^\infty{G_{\mu,\nu}((rr')^{\frac{1}{c}})f(r')\td r'}
 \end{equation*}
 defines a unitary involutive operator $T:L^2(\mathbb{R}_+)\longrightarrow L^2(\mathbb{R}_+)$ for $c=\frac{\mu+\nu+2}{2}>0$. Then the coordinate change $r=\left(\frac{x}{2}\right)^{2c}$, $r'=\left(\frac{y}{2}\right)^{2c}$ gives the claim.
 \item[\textup{(2)}] A short calculation using that $\mathcal{D}_{\mu,\nu}$ is a symmetric operator in $L^2(\mathbb{R}_+,x^{\mu+\nu+1}\td x)$ gives the desired statement if one knows that the $G$-function satisfies the following differential equation
 \begin{equation*}
  \left(\mathcal{D}_{\mu,\nu}\right)_xG_{\mu,\nu}\left(\left(\frac{xy}{4}\right)^2\right) = \left(\mathcal{D}_{\mu,\nu}\right)_yG_{\mu,\nu}\left(\left(\frac{xy}{4}\right)^2\right).
 \end{equation*}
 But this is easily derived from
 \begin{equation*}
  \theta_xG_{\mu,\nu}\left(\left(\frac{xy}{4}\right)^2\right) = 2\left(\theta G_{\mu,\nu}\right)\left(\left(\frac{xy}{4}\right)^2\right) = \theta_yG_{\mu,\nu}\left(\left(\frac{xy}{4}\right)^2\right)
 \end{equation*}
 using the expression \eqref{eq:DiffOp2} for $\mathcal{D}_{\mu,\nu}$ and the differential equation \eqref{eq:GDiffEq} for the $G$-function.\qedhere
\end{enumerate}
\end{pf}

\begin{theorem}[Meijer's $G$-transform]\label{thm:EigFctGTrafo}
Suppose that $\mu$ and $\nu$ satisfy \eqref{IntCond}. Then for each $j\in\mathbb{N}_0$ the function $\Lambda_{2,j}^{\mu,\nu}$ is an eigenfunction of Meijer's $G$-transform $\mathcal{T}_{\mu,\nu}$ with eigenvalue $(-1)^j$.
\end{theorem}

\begin{remark}
In Remark \ref{rem:DiffEqSpecialCase} we noted that for $\nu=-1$ the functions $\Lambda_{2,j}^{\mu,-1}$ are eigenfunctions of the self-adjoint second order differential operator $\mathcal{S}_{\mu,-1}$\index{Smupm1@$\mathcal{S}_{\mu,\pm1}$} on $L^2(\mathbb{R}_+,x^{\mu+\nu+1}\td x)$. In \cite[Theorem 4.1.1]{KM07a} Kobayashi and Mano proved that, if $\mu$ is an integer, the operator $e^{\frac{\pi\sqrt{-1}}{2}\mathcal{S}_{\mu,-1}}$, defined via spectral theory, is given as
\begin{equation*}
 e^{\frac{\pi\sqrt{-1}}{2}\mathcal{S}_{\mu,-1}} = e^{-\frac{\mu+1}{2}\pi\sqrt{-1}}\mathcal{T}_{\mu,-1}.
\end{equation*}
\end{remark}

One can use Theorem \ref{thm:EigFctGTrafo} to obtain an integral formula for the generating function $G_2^{\mu,\nu}$ and hence for the Bessel functions involved similar to the one in Corollary \ref{cor:BesselIntFormula}:

\begin{corollary}\label{cor:BesselIntFormula2}
Let $\frac{1}{2}<\alpha<\infty$, $\beta=\frac{\alpha}{2\alpha-1}$ and assume that $\mu$ and $\nu$ satisfy condition \eqref{IntCond}. Then for $x>0$
\begin{multline*}
 \int_0^\infty{G_{\mu,\nu}\left(\left(\frac{xy}{4}\right)^2\right)\widetilde{I}_{\frac{\mu}{2}}\left((\alpha-1)y\right)\widetilde{K}_{\frac{\nu}{2}}\left(\alpha y\right)y^{\mu+\nu+1}\td y}\\
 = 2^{\mu+\nu-1}\left(\frac{\beta}{\alpha}\right)^{\frac{\mu+\nu+2}{2}}\widetilde{I}_{\frac{\mu}{2}}((\beta-1)x)\widetilde{K}_{\frac{\nu}{2}}(\beta x).
\end{multline*}
\end{corollary}

\begin{pf}
By Theorem \ref{thm:EigFctGTrafo} we have
\begin{align*}
 \mathcal{T}_{\mu,\nu}\Lambda_{2,j}^{\mu,\nu} &= (-1)^j\Lambda_{2,j}^{\mu,\nu} \quad \text{ for every }j\in\mathbb{N}_0.
\end{align*}
Taking generating functions of both sides yields
\begin{equation*}
 \left(\mathcal{T}_{\mu,\nu}\right)_xG_2^{\mu,\nu}(t,x) = G_2^{\mu,\nu}(-t,x),
\end{equation*}
and this gives the desired formula for $\alpha=\frac{1}{1-t}$.
\end{pf}

\section{Applications of minimal representations}\label{sec:RepTh}

In this section we complete the proofs of Proposition \ref{prop:DiffOpProperties} (4) and Theorem \ref{thm:EigFctGTrafo} using representation theory.

There is a long history of research on the interactions of the theory of special functions and representation theory, in particular, related to the analysis on homogeneous spaces of Lie groups. Our fourth order differential operator $\mathcal{D}_{\mu,\nu}$ is also motivated by group representations. The general features in our setting, however, are quite different from the traditional analysis on homogeneous spaces.

In fact, we are dealing with an \lq extremely small\rq\ infinite dimensional representation, sometimes referred to as the {\textit{minimal representation}}, for which the whole group cannot act on the geometry of the $L^2$-model (the Gelfand--Kirillov dimension of the minimal representation is too small for that). Hence Lie algebra action here is not just given by vector fields, but involves higher order differential operators. (In our particular example second order suffices.) We shall prove in Theorem \ref{lem:CasimirAction} that the fourth order differential operator $\mathcal{D}_{\mu,\nu}$ comes from the Casimir operator when both $\mu$ and $\nu$ are integers $\geq-1$. We shall also discuss the underlying algebraic structure of various formulas on $(\Lambda_{2,j}^{\mu,\nu})_j$ such as recurrence relations (Theorem E) and Meijer's $G$-transforms (Theorem D) in connection with representation theory.

In the converse direction, the theory of special functions developed in this paper allows us to interpret the functions $\Lambda_{2,j}^{\mu,\nu}$ as a basis for the space $L^2(C)_{\textup{rad}}$ of $K'$-spherical vectors (see \eqref{eq:IsoL2rad}) and hence contributes to representation theory.

We adopt the notation of \cite{KM07b}. In particular, we assume that $p\geq q\geq2$ are integers of the same parity and $p+q\geq6$. With $\mu=p-3$, $\nu=q-3$ this is equivalent to $\mu\geq\nu\geq-1$ being integers with the same parity not both equal to $-1$, i.e. to the integrality condition \eqref{IntCond}.

Consider the indefinite orthogonal group \index{G@$G$|textbf}\index{O(p,q)@$O(p,q)$|textbf}
\begin{equation*}
 G:=O(p,q)=\{g\in\textup{GL}(p+q,\mathbb{R}):\ ^tgI_{p,q}g=I_{p,q}\}
\end{equation*}
of signature $(p,q)$, where \index{Ipq@$I_{p,q}$|textbf}
\begin{equation*}
 I_{p,q} := \left(\begin{array}{cc}I_p&0\\0&-I_q\end{array}\right).
\end{equation*}
Then $G$ is a semisimple non-compact Lie group. In \cite{KO03} Kobayashi and \O rsted construct an $L^2$-model of the minimal representation $\pi$\index{$\pi$|textbf} of $G$, which is an analogon of the classical Schr\"odinger model of the Shale--Segal--Weil representation of the metaplectic group. More precisely, it is an irreducible unitary representation on the Hilbert space $L^2(C)=L^2(C,\td\mu)$\index{L2C@$L^2(C)$|textbf}, where \index{C@$C$|textbf}
\begin{equation*}
 C = \{\zeta\in\mathbb{R}^{p+q-2}\ \backslash\ \{0\}:\zeta_1^2+\cdots+\zeta_{p-1}^2-\zeta_p^2-\cdots-\zeta_{p+q-2}^2=0\}
\end{equation*}
is an isotropic cone in $\mathbb{R}^{p+q-2}$ and $\td\mu$\index{dmu@$\td\mu$|textbf} an $O(p-1,q-1)$-invariant measure on $C$. Let $w_0:=I_{p,q}$\index{w0@$w_0$|textbf}. The unitary operator $\mathcal{F}_C:=\pi(w_0)$\index{FC@$\mathcal{F}_C$|textbf} plays a role of the Fourier transform on the isotropic cone $C$. We call $\mathcal{F}_C$ the \textit{unitary inversion operator}.

On the isotropic cone $C$, we set
\begin{equation*}
 r:=\sqrt{\zeta_1^2+\cdots+\zeta_{p-1}^2}=\sqrt{\zeta_p^2+\cdots\zeta_{p+q-2}^2}.
\end{equation*}
Then the pull-back of the projection $r:C\longrightarrow\mathbb{R}_+$ induces an isomorphism of Hilbert
spaces
\begin{equation}
 L^2(\mathbb{R}_+,\frac{1}{2}r^{p+q-5}\td r)\stackrel{\sim}{\longrightarrow}L^2(C)_{\textup{rad}}\subseteq L^2(C),\label{eq:IsoL2rad}
\end{equation}
where $L^2(C)_{\textup{rad}}$\index{L2Crad@$L^2(C)_{\textup{rad}}$|textbf} stands for the closed subspace of $L^2(C)$ consisting of radial functions, or equivalently, functions that are invariant by $K':=O(p-1)\times O(q-1)$\index{K'@$K'$|textbf}. Elements of $L^2(C)_{\textup{rad}}$ are called \textit{$K'$-spherical vectors}. To justify this name consider the bipolar coordinates which parameterize the cone $C$:
\begin{equation*}
 \mathbb{R}_+\times\mathbb{S}^{p-2}\times\mathbb{S}^{q-2}\longrightarrow C, (r,\omega,\eta)\mapsto(r\omega,r\eta),
\end{equation*}
where $\mathbb{S}^{n-1}$\index{Sn1@$\mathbb{S}^{n-1}$|textbf} denotes the unit sphere in $\mathbb{R}^n$. Then the action of $K'=O(p-1)\times O(q-1)$ is given by rotation in the spherical variables $\omega$ and $\eta$. Hence, $K'$-spherical functions are those depending only on the radial part $r$.

We take a negative definite bilinear form $B$\index{B@$B$|textbf} on $\mathfrak{o}(n)$ which is invariant under the adjoint action and turns $(E_{i,j}-E_{j,i})_{1\leq i<j\leq n}$ into an orthonormal basis with respect to $-B$. Restriction to the Lie algebra $\mathfrak{k}=\mathfrak{o}(p)\oplus\mathfrak{o}(q)$ of the maximal compact subgroup $K:=O(p)\times O(q)$\index{K@$K$|textbf} of $G$ yields the Casimir element\index{Omega@$\Omega$|textbf}
\begin{equation*}
 \Omega = -\sum_{0\leq i<j\leq p-1}{(E_{i,j}-E_{j,i})^2} -\sum_{p\leq i<j\leq p+q-1}{(E_{i,j}-E_{j,i})^2}.
\end{equation*}
The fourth order differential operator $\mathcal{D}_{\mu,\nu}$ with integral $\mu$ and $\nu$ arises from the Casimir action on the minimal representation as follows.

\begin{theorem}\label{lem:CasimirAction}
Suppose $p,q\geq2$. Under the coordinate change $x:=2r$ the operator $2\td\pi(\Omega)$ acts on $L^2(C)_{\textup{rad}}$ as the differential operator $\mathcal{D}_{p-3,q-3}$\index{Dmunu@$\mathcal{D}_{\mu,\nu}$}.
\end{theorem}

\begin{pf}
The group $K'$ acts trivially on $L^2(C)_{\textup{rad}}$ and so does its Lie algebra $\mathfrak{k}'$. Let $\Omega'$ be the Casimir element for $\mathfrak{k}'$. Then, the action of $\Omega$ on $L^2(C)_{\textup{rad}}$ is given by
\begin{align*}
 \td\pi(\Omega) &= \td\pi(\Omega) - \td\pi(\Omega')\\
 &= -\sum_{j=1}^{p-1}{(E_{0,j}-E_{0,j})^2}-\sum_{j=p}^{p+q-2}{(E_{j,p+q-1}-E_{p+q-1,j})^2}.
\end{align*}
By \cite[(2.3.15) and (2.3.19)]{KM07b} the Lie algebra action is given by differential operators in the ambient space $\mathbb{R}^{p+q-2}$ along $C$ as
\begin{align*}
 \td\pi(E_{0,j}-E_{j,0}) &= -\frac{\sqrt{-1}}{4}\left(4\zeta_j+\zeta_j\square-(2E+n-2)\frac{\partial}{\partial\zeta_j}\right)\\
 \intertext{for $1\leq j\leq p-1$ and similarly}
 \td\pi(E_{j,p+q-1}-E_{p+q-1,j}) &= -\frac{\sqrt{-1}}{4}\left(4\zeta_j+\zeta_j\square+(2E+n-2)\frac{\partial}{\partial\zeta_j}\right)
\end{align*}
for $p\leq j\leq p+q-2$, where \index{$\square$|textbf}\index{E@$E$|textbf}
\begin{align*}
 \square &= \frac{\partial^2}{\partial\zeta_1^2}+\cdots+\frac{\partial^2}{\partial\zeta_{p-1}^2}-\frac{\partial^2}{\partial\zeta_p^2}-\cdots-\frac{\partial^2}{\partial\zeta_{p+q-2}^2},\\
 E &= \sum_{j=1}^{p+q-2}{\zeta_j\frac{\partial}{\partial\zeta_j}}.
\end{align*}
Applying the operator $\td\pi(\Omega)$ to a radial function $f(r)$ gives the desired result.
\end{pf}

The reason why the Casimir operator $\td\pi(\Omega)$ is a differential operator of fourth order rather than of second order is the fact that the group $K=O(p)\times O(q)$ cannot act continuously in a non-trivial fashion on the cone $C$, whereas the smaller group $K'\subseteq K$ does.

If one restricts the representation to the maximal compact subgroup $K=O(p)\times O(q)\subseteq G$\index{K@$K$} one has the following multiplicity-free decomposition of $\pi$ into $K$-irreducible representations
\begin{equation}
 L^2(C) \cong {\sum_{a=0}^\infty}\raisebox{0.15cm}{$^\oplus$}{\,V^a}\label{eq:KtypeDecomp}
\end{equation}
with \index{Va@$V^a$|textbf}
\begin{equation*}
 V^a \cong \mathcal{H}^a(\mathbb{R}^p)\otimes\mathcal{H}^{a+\frac{p-q}{2}}(\mathbb{R}^q).
\end{equation*}
Here $\mathcal{H}^j(\mathbb{R}^n)$\index{Hj(Rn)@$\mathcal{H}^j(\mathbb{R}^n)$|textbf} is the space of spherical harmonics on $\mathbb{R}^n$ of degree $j$, on which the orthogonal group $O(n)$ acts irreducibly. The branching laws for $O(n)\downarrow O(n-1)$ ($n\geq2$) show that the $O(n-1)$-fixed vectors in $\mathcal{H}^j(\mathbb{R}^n)$ form a one-dimensional subspace. Hence, if one restricts to the smaller subgroup $K'=O(p-1)\times O(q-1)\subseteq K$\index{K'@$K'$}, then the $K'$-spherical vectors in each $K$-type $V^a$ form a one-dimensional subspace $V^a\cap L^2(C)_{\textup{rad}}$. Considering only the radial functions in \eqref{eq:KtypeDecomp} we obtain
\begin{equation}
 L^2(C)_{\textup{rad}} \cong {\sum_{a=0}^\infty}\raisebox{0.15cm}{$^\oplus$}{\,\left(V^a\cap L^2(C)_{\textup{rad}}\right)},\label{eq:KtypeDecompRad}
\end{equation}
because the property of being radial is determined by the subgroup $K'$ of $K$. We will show that each subspace $V^a\cap L^2(C)_{\textup{rad}}$ is spanned by $\Lambda_{2,a}^{p-3,q-3}(2r)$.

In order to do so we observe that by Schur's Lemma the Casimir element $\Omega$ acts as a scalar on each $K$-type $V^a$, because $\Omega$ is central in the universal enveloping algebra $\mathcal{U}(\mathfrak{k})$. In fact, representation theory of the compact group $K=O(p)\times O(q)$ shows that $\Omega$ acts on $V^a$ by
\begin{equation}
 a(a+p-2)+\left(a+\frac{p-q}{2}\right)\left(a+\frac{p+q-4}{2}\right).\label{eq:CasimirEigVal}
\end{equation}

Now we are ready to prove that the differential operator $\mathcal{D}_{\mu,\nu}$ is self-adjoint.

\begin{proof}[\sc Proof of Proposition \ref{prop:DiffOpProperties}~(4)]
By \cite[Theorem 4.4.4.3]{War72} the Casimir element $\Omega$ defines a self-adjoint operator on the representation space $L^2(C)$. On the other hand, the element $-2\Omega$ acts on $L^2(C)_{\textup{rad}}$ as $\mathcal{D}_{\mu,\nu}$ by Theorem \ref{lem:CasimirAction}. Hence $\mathcal{D}_{\mu,\nu}$ extends to a self-adjoint operator on the space $L^2(C)_{\textup{rad}}\cong L^2(\mathbb{R}_+,x^{\mu+\nu+1}\td x)$ of radial functions. The decomposition \eqref{eq:KtypeDecompRad} is discrete, so the spectrum of $\mathcal{D}_{\mu,\nu}$ is discrete and given by twice the eigenvalues \eqref{eq:CasimirEigVal} of the Casimir element which in view of $j=a, \mu=p-3$, and $\nu=q-3$ are given by \index{$\lambda_j^{\mu,\nu}$}
\begin{equation*}
 \lambda_j^{\mu,\nu} = 2j(j+\mu+1) + 2\left(j+\frac{\mu-\nu}{2}\right)\left(j+\frac{\mu+\nu+2}{2}\right).
\end{equation*}
Since the spaces $V^a\cap L^2(C)_{\textup{rad}}$ in \eqref{eq:KtypeDecompRad} are one-dimensional and the $(\lambda_j^{\mu,\nu})_{j\in\mathbb{N}_0}$ is strictly increasing, the $L^2$-eigenspaces of $\mathcal{D}_{\mu,\nu}$ are also one-dimensional.
\end{proof}

For the $L^2$-model of the minimal representation, it is in general a difficult problem to find explicit $K$-finite vectors (see \cite[Introduction]{KM07a}). Theorem \ref{lem:CasimirAction} implies that the functions $\Lambda_{2,j}^{\mu,\nu}$ ($j=0,1,2,\ldots$) for integers $\mu\geq\nu$ of the same parity are indeed $K$-finite vectors of the minimal representation $\pi$. More precisely, we have:

\begin{corollary}\label{cor:Kfin}
For every $j\in\mathbb{N}_0$, the space $V^j\cap L^2(C)_{\textup{rad}}$ is one-dimensional. It is spanned by $\Lambda_{2,j}^{p-3,q-3}(2r)$.
\end{corollary}

\begin{pf}
We have already seen in Theorem \ref{thm:EigFct} that the functions $\Lambda_{2,j}^{p-3,q-3}(x)$ are eigenfunctions of $\mathcal{D}_{p-3,q-3}$ in $L^2(\mathbb{R}_+,\frac{1}{2}r^{p+q-5}\td r)$ for the eigenvalues $\lambda_j^{p-3,q-3}$. Hence, $\Lambda_{2,j}^{p-3,q-3}(2r)$ is an eigenfunction of the Casimir operator $\td\pi(\Omega)$ for the eigenvalue \eqref{eq:CasimirEigVal}. Since all those eigenvalues are distinct, the claim follows.
\end{pf}

\begin{remark}
\begin{enumerate}
 \item[\textup{(1)}] For $q=2$ ($p$ and $j$ arbitrary), we saw in Corollary \ref{cor:SpecialValue} that
 \begin{equation*}
   \Lambda_{2,j}^{p-3,q-3}(2r) = \const\cdot L_j^{p-3}(4r)e^{-2r}.
 \end{equation*}
 In this special case Corollary \ref{cor:Kfin} was already proved in \cite[Proposition 3.2.1~(2)]{KM07a}.
 \item[\textup{(2)}] For $j=0$ ($p$ and $q$ arbitrary) Example \ref{ex:3Fcts} gives
 \begin{equation*}
  \Lambda_{2,0}^{p-3,q-3}(2r) = \const\cdot \widetilde{K}_{\frac{q-3}{2}}(2r).
 \end{equation*}
 This bottom case of Corollary \ref{cor:Kfin} was already found in \cite[Theorem 5.8~(1)]{KO03} by taking theFourier transform of the \lq conformal model\rq\ of the minimal representation.
\end{enumerate}
\end{remark}

The relation of the $G$-transform considered in Section \ref{sec:GTrafo} to representation theory was found by Kobayashi and Mano in \cite{KM07bprocj,KM07b}.

\begin{theorem}[{\cite[Theorem 4.1.1]{KM07b}}]\label{thm:WeilInv}
Under the coordinate change $x=2r$ the unitary involution operator $\mathcal{F}_C$\index{FC@$\mathcal{F}_C$} acts on $L^2(C)_{\textup{rad}}$ as the $G$-transform $(-1)^{\frac{p-q}{2}}\mathcal{T}_{p-3,q-3}$.\index{Tmunu@$\mathcal{T}_{\mu,\nu}$}
\end{theorem}

In the case that $\mu=p-3$ and $\nu=q-3$ satisfy \eqref{IntCond} we can now give an alternative proof of the results of Proposition \ref{prop:GTrafoProperties} using representation theory. By Theorem \ref{thm:WeilInv} the operator $\mathcal{T}_{p-3,q-3}$ is given as the restriction of the unitary inversion operator $\mathcal{F}_C=\pi(w_0)$ to $L^2(C)_{\textup{rad}}$. This shows that $\mathcal{T}_{p-3,q-3}$ defines a unitary operator on $L^2(\mathbb{R}_+,x^{\mu+\nu+1}\td x)$. Furthermore, since $w_0$ normalizes (actually, centralizes) the Lie algebra $\mathfrak{k}$, $w_0$ commutes with the Casimir operator $\Omega$. Therefore, in view of Theorems \ref{lem:CasimirAction} and \ref{thm:WeilInv}, also $\mathcal{D}_{p-3,q-3}$ and $\mathcal{T}_{p-3,q-3}$ commute.

Theorem \ref{thm:WeilInv} also allows us to give the announced proof of the fact that the functions $\Lambda_{2,j}^{p-3,q-3}$ are eigenfunctions of Meijer's $G$-transform. 
\begin{proof}[\sc Proof of Theorem \ref{thm:EigFctGTrafo}]
It is observed in \cite[Proof of Theorem 3.1.1~(3)]{KM07b} that $w_0$ acts on each $K$-type $V^a$ as $(-1)^{a+\frac{p-q}{2}}$. Since $\Lambda_{2,a}^{p-3,q-3}(2r)$ is a $K$-finite vector in the $K$-type $V^a$ by Corollary \ref{cor:Kfin}, the statement follows from Theorem \ref{thm:WeilInv}.
\end{proof}

Finally, we can also give a representation theoretic explanation for the recurrence relations in Propositions \ref{prop:RecRelH} and \ref{prop:RecRelxSq}. Using the notation of \cite{KM07b} we consider the Langlands decomposition of the maximal parabolic subalgebra $\overline{\mathfrak{p}^{\textup{max}}}$ of $\mathfrak{g}=\Lie(G)$\index{g@$\mathfrak{g}$|textbf}:
\begin{equation*}
 \overline{\mathfrak{p}^{\textup{max}}} = \mathfrak{m}^{\textup{max}}+\mathfrak{a}+\overline{\mathfrak{n}^{\textup{max}}}.
\end{equation*}
$\mathfrak{m}^{\textup{max}}$ is the Lie algebra of $O(p-1,q-1)$ which acts on the geometry $C$. The abelian part $\mathfrak{a}$ is one-dimensional and spanned by $H\in\mathfrak{a}$ which acts on $L^2(C)$ as the first order differential operator
\begin{equation*}
 \td\pi(H)=-\left(E+\frac{p+q-4}{2}\right),
\end{equation*}
where $E$\index{E@$E$} is again the Euler operator. Hence $H$ acts on radial functions as the differential operator $-\mathcal{H}_{p+q-6}$\index{Halpha@$\mathcal{H}_\alpha$} (see Section \ref{sec:RecRel}). Finally $\overline{\mathfrak{n}^{\textup{max}}}$ is spanned by $\overline{N}_1,\ldots,\overline{N}_{p+q-2}\in\overline{\mathfrak{n}^{\textup{max}}}$ whose action is given by the multiplication operators
\begin{equation*}
 \td\pi(\overline{N}_j) = 2\sqrt{-1}\zeta_j.
\end{equation*}
In particular, the sum of squares $\sum_{j=1}^{p+q-2}{\td\pi(\overline{N}_j)^2}=-8r^2$ leaves the space of radial functions $L^2(C)_{\textup{rad}}$ invariant.

The key to an understanding of the underlying algebraic structure of the recurrence relations in Propositions \ref{prop:RecRelH} and \ref{prop:RecRelxSq} is the action of $H$ and $\overline{N}_j$ on the $K$-types $V^a$. For convenience put $V^{-1}:=0$.

\begin{lemma}
The Lie algebra action $\td\pi(X):\bigoplus_{a=0}^\infty{V^a}\longrightarrow\bigoplus_{a=0}^\infty{V^a}$ ($X\in\mathfrak{g}$) induces the following linear maps for each $a\in\mathbb{N}_0$:
\begin{align*}
 \td\pi(H):V^a &\longrightarrow V^{a+1}\oplus V^{a-1},\\
 \td\pi(\overline{N}_j):V^a &\longrightarrow V^{a+1}\oplus V^a\oplus V^{a-1}, & 1\leq j\leq p+q-2.
\end{align*}
\end{lemma}

\begin{pf}
Let $\mathfrak{g}=\mathfrak{k}+\mathfrak{p}$ be the Cartan decomposition of $\mathfrak{g}$ determined by the choice of $\mathfrak{k}$. For $X\in\mathfrak{k}$ we have
\begin{equation}
 \td\pi(X):V^a\longrightarrow V^a\label{eq:MapPropK}
\end{equation}
since $V^a$ is a $K$-module. For the action of $\mathfrak{p}$ consider the $K$-module homomorphism
\begin{equation*}
 \mathfrak{p}\otimes V^a\rightarrow L^2(C),\ X\otimes\phi\mapsto\td\pi(X)\phi,
\end{equation*}
where $K$ acts on $\mathfrak{p}$ by the adjoint action. To determine which $K$-types appear in the image one has to decompose $\mathfrak{p}\otimes V^a$ into simple $K$-modules. Let $\boxtimes$ denote the outer tensor product representation. Since $\mathfrak{p}\cong\mathcal{H}^1(\mathbb{R}^p)\boxtimes\mathcal{H}^1(\mathbb{R}^q)$ as $K$-modules, we find
\begin{align*}
 \mathfrak{p}\otimes V^a &\cong (\mathcal{H}^1(\mathbb{R}^p)\otimes\mathcal{H}^a(\mathbb{R}^p))\boxtimes(\mathcal{H}^1(\mathbb{R}^q)\otimes\mathcal{H}^{a+\frac{p-q}{2}}(\mathbb{R}^q)).
\end{align*}
Using branching rules for tensor products of $O(n)$-modules (see e.g. \cite{HTW05}) one obtains that
\begin{equation*}
 \mathcal{H}^1(\mathbb{R}^n)\otimes\mathcal{H}^a(\mathbb{R}^n) \cong \mathcal{H}^{a+1}(\mathbb{R}^n)\oplus\mathcal{H}^{a-1}(\mathbb{R}^n)\oplus W,
\end{equation*}
where $\mathcal{H}^{-1}(\mathbb{R}^n):=0$ and $W$ is either $0$ or an irreducible $O(n)$-module which is not isomorphic to any of the $\mathcal{H}^k(\mathbb{R}^n)$. Hence
\begin{align*}
 \mathfrak{p}\otimes V^a &\cong (\mathcal{H}^{a+1}(\mathbb{R}^p)\oplus\mathcal{H}^{a-1}(\mathbb{R}^p)\oplus W_1)\boxtimes(\mathcal{H}^{a+1+\frac{p-q}{2}}(\mathbb{R}^q)\oplus\mathcal{H}^{a-1+\frac{p-q}{2}}(\mathbb{R}^q)\oplus W_2)\\
 &\cong (\mathcal{H}^{a+1}(\mathbb{R}^p)\boxtimes\mathcal{H}^{a+1+\frac{p-q}{2}}(\mathbb{R}^q))\oplus(\mathcal{H}^{a-1}(\mathbb{R}^p)\boxtimes\mathcal{H}^{a-1+\frac{p-q}{2}}(\mathbb{R}^q))\oplus W'\\
 &\cong V^{a+1}\oplus V^{a-1}\oplus W',
\end{align*}
where $W'$ is a sum of simple $K$-modules that do not appear in $\pi$, i.e. different from the $V^j$, $j\in\mathbb{N}_0$. Hence for $X\in\mathfrak{p}$ we have
\begin{equation}
 \td\pi(X):V^a\longrightarrow V^{a+1}\oplus V^{a-1}.\label{eq:MapPropP}
\end{equation}
Putting \eqref{eq:MapPropK} and \eqref{eq:MapPropP} together proves the claim since $H\in\mathfrak{a}\subseteq\mathfrak{p}$.
\end{pf}

Note that both $\td\pi(H)$ and $\sum_{j=1}^{p+q-2}{\td\pi(\overline{N}_j)^2}$ leave $L^2(C)_{\textup{rad}}$ invariant. Since $V^a\cap L^2(C)_{\textup{rad}}$ is one-dimensional and spanned by $\Lambda_{2,a}^{p-3,q-3}(2r)$ for every $a\in\mathbb{N}_0$, we obtain
\begin{align*}
 \mathcal{H}_{p+q-6}\Lambda_{2,j}^{p-3,q-3} &\in\textup{span}\{\Lambda_{2,k}^{p-3,q-3}: k=j-1,j+1\},\\
 x^2\Lambda_{2,j}^{p-3,q-3} &\in\textup{span}\{\Lambda_{2,k}^{p-3,q-3}:k=j-2,j-1,j,j+1,j+2\}.
\end{align*}
which can be viewed as qualitative versions of Propositions \ref{prop:RecRelH} and \ref{prop:RecRelxSq}.

\addcontentsline{toc}{section}{References}
\providecommand{\bysame}{\leavevmode\hbox to3em{\hrulefill}\thinspace}
\providecommand{\MR}{\relax\ifhmode\unskip\space\fi MR }
\providecommand{\MRhref}[2]{\href{http://www.ams.org/mathscinet-getitem?mr=#1}{#2}}
\providecommand{\href}[2]{#2}

\vspace{30pt}

\textsc{Joachim Hilgert\\Institut f\"ur Mathematik, Universit\"at Paderborn, Warburger Str. 100, 33098 Paderborn, Germany.}\\
\textit{E-mail address:} \texttt{hilgert@math.uni-paderborn.de}\\

\textsc{Toshiyuki Kobayashi}\\
\textit{Home address:} \textsc{Graduate School of Mathematical Sciences, IPMU, the University of Tokyo, 3-8-1 Komaba, Meguro, Tokyo, 153-8914, Japan.}\\
\textit{Current address:} \textsc{Max-Planck-Institut f\"ur Mathematik, Vivatsgasse 7, 53111 Bonn, Germany.}\\
\textit{E-mail address:} \texttt{toshi@ms.u-tokyo.ac.jp}\\

\textsc{Gen Mano\\Graduate School of Mathematical Sciences, the University of Tokyo, 3-8-1 Komaba, Meguro, Tokyo, 153-8914, Japan.}\\

\textsc{Jan M\"ollers\\Institut f\"ur Mathematik, Universit\"at Paderborn, Warburger Str. 100, 33098 Paderborn, Germany.}\\
\textit{E-mail address:} \texttt{moellers@math.uni-paderborn.de}

\addcontentsline{toc}{section}{Index}

\begin{theindex}

  \item $(a)_n$, \textbf{11}
  \item $\Lambda_{i,j}^{\mu,\nu}(x)$, \hyperpage{3}, \textbf{10}
  \item $\delta(i)$, \hyperpage{4}, \textbf{9}, \hyperpage{22}
  \item $\epsilon(i)$, \hyperpage{4}, \textbf{9}, \hyperpage{22}
  \item $\lambda_j^{\mu,\nu}$, \textbf{3}, \hyperpage{6}, 
		\hyperpage{26}
  \item $\pi$, \textbf{24}
  \item $\square$, \textbf{25}
  \item $\theta$, \textbf{2}

  \indexspace

  \item $a_\alpha$, \textbf{9}, \hyperpage{14}

  \indexspace

  \item $B$, \textbf{25}
  \item $b_\alpha$, \textbf{9}, \hyperpage{14}

  \indexspace

  \item $C$, \hyperpage{5}, \textbf{24}
  \item $c_{1,j}^{\mu,\nu}$, \hyperpage{4}, \textbf{15}
  \item $c_{2,j}^{\mu,\nu}$, \hyperpage{4}, \textbf{15}
  \item $C_{\mu,\pm1}$, \textbf{6}, \hyperpage{18}

  \indexspace

  \item $\td\mu$, \textbf{24}
  \item $\mathcal{D}_{\mu,\nu}$, \textbf{2}, \hyperpage{6}, 
		\hyperpage{25}

  \indexspace

  \item $E$, \textbf{25}, \hyperpage{28}

  \indexspace

  \item ${_2F_2}(a_1,a_2;b_1,b_2;z)$, \hyperpage{18}
  \item $\mathcal{F}_C$, \textbf{24}, \hyperpage{27}

  \indexspace

  \item $G$, \hyperpage{5}, \textbf{24}
  \item $\mathfrak{g}$, \textbf{27}
  \item $G^{20}_{04}(t\vert b_1,b_2,b_3,b_4)$, \hyperpage{4}, 
		\textbf{22}
  \item $G_i^{\mu,\nu}(t,x)$, \hyperpage{2}, \textbf{7}
  \item $G_{\mu,\nu}(t)$, \hyperpage{4}, \textbf{22}

  \indexspace

  \item $\mathcal{H}_\alpha$, \textbf{19}, \hyperpage{28}
  \item $\mathcal{H}^j(\mathbb{R}^n)$, \textbf{26}
  \item $\mathcal{H}_\pm^{\mu,\nu}$, \textbf{4}

  \indexspace

  \item $\widetilde{I}_\alpha(z)$, \hyperpage{2}, \textbf{7}, 
		\hyperpage{16}
  \item $I_{p,q}$, \textbf{24}

  \indexspace

  \item $J_\alpha(z)$, \textbf{22}

  \indexspace

  \item $K$, \textbf{25}, \hyperpage{26}
  \item $K'$, \textbf{25}, \hyperpage{26}
  \item $\widetilde{K}_\alpha(z)$, \hyperpage{2}, \textbf{7}, 
		\hyperpage{16}

  \indexspace

  \item $L^2(C)$, \hyperpage{5}, \textbf{24}
  \item $L^2(C)_{\textup{rad}}$, \textbf{25}
  \item $L_n^\alpha(x)$, \hyperpage{3}, \textbf{15}, \hyperpage{17}

  \indexspace

  \item $O(p,q)$, \hyperpage{5}, \textbf{24}
  \item $\Omega$, \textbf{25}

  \indexspace

  \item $\mathcal{S}_{\mu,\pm1}$, \textbf{6}, \hyperpage{18}, 
		\hyperpage{23}
  \item $\mathbb{S}^{n-1}$, \hyperpage{5}, \textbf{25}

  \indexspace

  \item $\mathcal{T}_{\mu,\nu}$, \hyperpage{4}, \textbf{22}, 
		\hyperpage{27}

  \indexspace

  \item $V^a$, \textbf{26}

  \indexspace

  \item $w_0$, \textbf{24}

\end{theindex}

\end{document}